\documentclass[pamm,a4paper,fleqn]{w-art}
\usepackage{times,cite,w-thm}
\usepackage{hyperref}
\usepackage[capitalize]{cleveref}
\hypersetup{
   colorlinks=true,
   urlcolor=green!70!black,
   citecolor=blue!70!black,
   linkcolor=green!70!black%
}
\usepackage[T1]{fontenc}
\usepackage[utf8]{inputenc}
\theoremstyle{plain}

\theoremstyle{definition}

\usepackage{graphicx}
\usepackage{xcolor}
\usepackage{booktabs}
\usepackage{tikz}
\usetikzlibrary{calc}
\usepackage{amsmath,amsthm,amssymb,mathtools}
\usepackage{colortbl}
\usepackage{array}

\newcolumntype{L}[1]{>{\raggedright\let\newline\\\arraybackslash\hspace{0pt}}m{#1}}
\newcolumntype{C}[1]{>{\centering\let\newline\\\arraybackslash\hspace{0pt}}m{#1}}
\newcolumntype{R}[1]{>{\raggedleft\let\newline\\\arraybackslash\hspace{0pt}}m{#1}}
\usepackage{nicematrix}

\newcommand{\VEM}{\mathbb{V}}
\newcommand{\TT}{\mathbb{T}}
\newcommand{\Th}{\mathcal{T}_h}

\newcommand{\Fhkz}{F_{h,0}^k}
\newcommand{\Fhkd}{F_{h,D}^k}
\newcommand{\Vhk}{V_h^k}

\newcommand{\Whkz}{W_{h,0}^k}
\newcommand{\Whkd}{W_{h,D}^k}

\newcommand{\avg}[1]{\left\{\!\!\left\{ #1 \right\}\!\!\right\}}
\newcommand{\jump}[1]{[\![ #1 ]\!]}
\newcommand{\ncdof}{\texttt{ncdof}\xspace}
\newcommand{\nnze}{\texttt{nnze}\xspace}
\newcommand{\ndof}{\texttt{ndof}\xspace}

\newcommand{\Nel}{N_{\text{El}}}
\newcommand{\Nft}{N_{\text{Ft}}}

\newcommand{\Nb}[2]{N^{\text{#2,nb}}_{#1}}

\definecolor{pscol}{rgb}{0,0.6,0}

\usepackage{tikz,tikz-3dplot}
\usetikzlibrary{backgrounds,3dtools}
\newsavebox\TruncatedOctahedronTeal
\newsavebox\TruncatedOctahedronOrange
\tdplotsetmaincoords{80}{105}
\newcommand{\TruncatedOctahedron}{%
\begin{tikzpicture}[tdplot_main_coords,line cap=round,line join=round]
\path foreach \Coord [count=\X] in 
{(-1.5,-0.5,0.), (-1.5,0.5,0.), (-1.,-1.,-0.707107), (-1.,-1.,0.707107), 
(-1.,1.,-0.707107), (-1.,1.,0.707107), (-0.5,-1.5,0.), (-0.5, -0.5,-1.41421),
(-0.5,-0.5,1.41421), (-0.5,0.5,-1.41421), (-0.5,0.5, 1.41421), (-0.5,1.5,0.),
(0.5,-1.5,0.), (0.5,-0.5,-1.41421), (0.5,-0.5, 1.41421), (0.5,0.5,-1.41421),
(0.5,0.5,1.41421), (0.5,1.5,0.), (1.,-1., -0.707107), (1.,-1.,0.707107),
(1.,1.,-0.707107), (1.,1.,0.707107), (1.5, -0.5,0.), (1.5,0.5,0.)}
  {\Coord coordinate (p\X) \pgfextra{\xdef\NumVertices{\X}}};
 \path[overlay] ({sin(\tdplotmaintheta)*sin(\tdplotmainphi)},
    {-1*sin(\tdplotmaintheta)*cos(\tdplotmainphi)},
    {cos(\tdplotmaintheta)}) coordinate (n)
    (0.5,0.5,{0.5*sqrt(2)})  coordinate (L); 
 \edef\lstPast{0}
 \foreach \poly in 
 {{17,11,9,15}, {14,8,10,16}, {22,24,21,18}, {12,5,2,6}, {13,19,23,20},
    {4,1,3,7}, {19,13,7,3, 8,14}, {15,9,4,7,13,20}, {16,10,5,12,18,21}, 
    {22,18,12,6,11,17}, {20,23,24,22,17,15}, {14,16,21,24,23, 19}, {9,11,6,2,1,4},  
    {3,1,2,5,10,8}}
 {
  \pgfmathtruncatemacro{\ione}{{\poly}[0]}
  \pgfmathtruncatemacro{\itwo}{{\poly}[1]}
  \pgfmathtruncatemacro{\ithree}{{\poly}[2]}
  \path[overlay,3d coordinate={(dA)=(p\itwo)-(p\ione)},
   3d coordinate={(dB)=(p\itwo)-(p\ithree)},
   3d coordinate={(nA)=(dA)x(dB)}] ;
  \pgfmathtruncatemacro{\jtest}{sign(TD("(nA)o(p\ione)"))}
  \ifnum\jtest<0
   \path[overlay,3d coordinate={(nA)=(dB)x(dA)}];
  \fi
  \pgfmathsetmacro\myproj{TD("(nA)o(n)")}
  \pgfmathsetmacro\lproj{TD("(nA)o(L)")}
  \pgfmathtruncatemacro{\itest}{sign(\myproj)}
  \pgfmathtruncatemacro{\cf}{70+20*\lproj}%
  \ifnum\itest>-1 
   \draw[ultra thin] [fill=mypolyhedroncolor!\cf]
   plot[samples at=\poly,variable=\x](p\x) -- cycle; 
  \else
   \begin{scope}[on background layer] 
    \draw[gray!20,ultra thin] [fill=mypolyhedroncolor!\cf!black]
    plot[samples at=\poly,variable=\x](p\x) -- cycle;  
   \end{scope}
  \fi
  }  
\end{tikzpicture}}
\colorlet{mypolyhedroncolor}{teal}
\sbox\TruncatedOctahedronTeal{\TruncatedOctahedron}
\colorlet{mypolyhedroncolor}{orange}
\sbox\TruncatedOctahedronOrange{\TruncatedOctahedron}
\DeclareRobustCommand\textcircle[1]{\tikz[baseline=-2.5pt]{
            \node[shape=circle,draw,inner sep=2pt,#1] (char) {};}}

\begin{document}

\TitleLanguage[EN]
\title[Sparsity comparison of polytopal finite element methods]{Sparsity comparison of polytopal finite element methods}

\author{\firstname{Christoph}  \lastname{Lehrenfeld}\inst{1,}%
  \footnote{Corresponding author: e-mail \ElectronicMail{lehrenfeld@math.uni-goettingen.de},
            phone +49\,551\,39\,24565,
            fax +49\,551\,39\,33944}}

\address[\inst{1}]{\CountryCode[DE]Institut für Numerische und Angewandte Mathematik, Georg-August Universität Göttingen, Lotzestraße 16-18, 37083 Göttingen}
\author{\firstname{Paul} \lastname{Stocker}\inst{2,}%
  \footnote{e-mail \ElectronicMail{paul.stocker@univie.ac.at}}}

\address[\inst{2}]{\CountryCode[AT]Institut für Mathematik, Universität Wien, Oskar-Morgenstern-Platz 1, 1090 Wien}
\author{\firstname{Maximilian} \lastname{Zienecker}\inst{1,}%
\footnote{e-mail \ElectronicMail{m.zienecker@stud.uni-goettingen.de}}}
\AbstractLanguage[EN]
\begin{abstract}
  In this work we compare crucial parameters for efficiency of different finite element methods for solving partial differential equations (PDEs) on polytopal meshes. 
  We consider the Virtual Element Method (VEM) and different Discontinuous Galerkin (DG) methods, namely the Hybrid DG and Trefftz DG methods.
  The VEM is a conforming method, that can be seen as a generalization of the classic finite element method to arbitrary polytopal meshes.
  DG methods are non-conforming methods that offer high flexibility, but also come with high computational costs.
  Hybridization reduces these costs by introducing additional facet variables, onto which the computational costs can be transfered to. Trefftz DG methods achieve a similar reduction in complexity by selecting a special and smaller set of basis functions on each element. 
  The association of computational costs to different geometrical entities (elements or facets) leads to differences in the performance of these methods on different grid types. This paper aims to compare the dependency of these approaches across different grid configurations.
\end{abstract}
\maketitle                   %

\section{Introduction}
Finite element methods on polytopal meshes have gained increasing attention in recent years due to their flexibility and efficiency in handling complex geometries. 
Allowing general polytopic elements provides enormous flexibility in mesh generation, for meshing complex geometries using a minimum number of elements, for mesh adaptation and mesh coarsening.
In this work, we compare finite element methods for polytopal meshes on different mesh topologies. 

The Discontinuous Galerkin (DG) methods, see e.g. \cite{arnold2002unified, di2011mathematical}, use non-conforming basis functions defined on each element separately. 
Consistency across element interfaces is enforced by numerical fluxes. 
This offers versatile solutions for a wide range of problems and allows for polytopal elements and hanging nodes in the mesh.
However, compared to conforming methods, DG methods often incur significantly higher computational costs. This is primarily due to the emergence of numerous additional degrees of freedom (dofs) and a notable increase in non-zero entries in linear equation systems.

To mitigate these computational expenses, the principle of hybridization has been introduced \cite{cockburn2009unified}. Hybridization involves introducing \emph{facet unknowns} on the mesh skeleton in addition to existing \emph{volume unknowns}, and then reducing all local volume unknowns to the facet unknowns through static condensation. This approach proves particularly effective for higher-order methods, as the dominant computational costs scale with the globally coupled degrees of freedom which itself scale with $\mathcal{O}(k^{d-1})$ instead of $\mathcal{O}(k^{d})$, where $k$ denotes the polynomial degree of discretization and $d$ signifies the spatial dimension.

An alternative approach, which achieves a similar reduction in complexity, is the concept of Trefftz DG methods, see e.g. \cite{TrefftzSurvey}. 
In this method, basis functions on each element are adaptively chosen to conform to the underlying PDE problem (without including boundary or element interface conditions), an idea originating from \cite{trefftz1926}. 
As a result, the same approximation order for solutions of the PDE is attained with fewer unknowns, also scaling with $\mathcal{O}(k^{d-1})$ instead of $\mathcal{O}(k^{d})$.
The method has recently been applied to a wide range of problems including inhomogeneous right-hand sides and varying coefficients, see e.g.  \cite{LS_IJMNE_2023,IGMS_MC_2021,LLS_NM_2024,perinati2023quasi}.%

The Virtual Element Method (VEM), see e.g. \cite{vembasicprinciples, vemhitchhiker, vembook}, is a conforming method that can be seen as a generalization of the classical finite element method -- which can only be applied to meshes consisting of a restricted set of element types -- to arbitrary polytopal meshes. 
The VEM makes it possible to construct discrete spaces featuring (higher-order) continuity properties on polytopal meshes, however, the basis functions are not explicitly defined but are only known through their degrees of freedom.

Besides the number of globally coupled unknowns, the coupling pattern is crucial for the performance of a numerical method. For the different approaches addressed before, the coupling pattern is different and varies across different grid types.
In this paper, we aim to compare the dependency of the different approaches across different grid types. By examining their performance under various grid configurations, we seek to elucidate their relative merits and drawbacks, providing insights into their suitability for different computational scenarios. 

The remainder of this paper is structured as follows: In \cref{sec:model_problem}, we briefly recap the discretization methods on the example of a simple Laplace problem. In \cref{sec:measures}, we introduce notation and measures that are to be investigated in \cref{sec:periodicpolys} for different grid types. We conclude with a discussion of the results in \cref{sec:conclusion}.

\section{Polytopal finite element methods for a model problem} \label{sec:model_problem}
To recap the essentials of the different discretizations, we consider the Laplace problem as a simple model problem and summarize the approaches for that example. 
For the sake of simplicity in the presentation, we consider the Laplace problem with inhomogeneous Dirichlet boundary conditions on the unit square, respectively unit cube, $\Omega \in \{ (0,1)^2, (0,1)^3 \}$ for the discussion of the methods here, while for the comparison of methods later on we will switch to periodic meshes and periodic boundary conditions to get rid of boundary effects.

\subsection{Elliptic model problem}
We consider the following model problem: 
Find $u: \Omega \to \mathbb{R}$ such that
\begin{equation}
  - \Delta u = 0 \quad \text{in } \Omega, \qquad
  \quad u = u_D \quad \text{on } \partial \Omega.
\end{equation}
The weak formulation reads: Find $u \in H^1_{D}(\Omega)$, such that
\begin{align}
  a(u,v) &= 0 \quad \forall v \in H^1_0(\Omega) \quad \text{ with }\quad
  a(u,v) = (\nabla u,\nabla v)_{\Omega}, %
\end{align}
where $(\cdot,\cdot)_S$ denotes the $L^2$-inner product over the set $S$ and $H^1_0(\Omega)$ and 
$H^1_D(\Omega)$ are the Sobolev spaces with zero boundary trace and Dirichlet boundary traces, respectively.
We assume that $\Omega$ is decomposed into a triangulation $\Th$ of the domain $\Omega$ consisting of $N_{\text{El}}$ elements and a set of  $N_{\text{Ft}}$ facets $\mathcal{F}_h$. The elements $T \in \Th$ are assumed to be (straight) polytopal, for instance (but not limited to) triangles, quadrilaterals, hexagons in 2D and for instance tetrahedra, or hexahedra in 3D.

The application of standard conforming finite element methods is not possible on general polytopal meshes. In the subsequent sections, we will consider different discretization methods that are non-conforming except for the Virtual Element Method (VEM). 

\begin{figure}
	\def\figheight{1.25cm}
	\def\figlength{1.5*\figheight}
	\def\noderadius{\figlength/1.5cm}
	\def\diff{1.5mm}
	\begin{minipage}[t]{.31\textwidth}
        \centering
		\begin{tikzpicture}
			\draw (-\diff,0) -- (-\diff,0.75*\figheight) -- (-\figlength/2-\diff,1.25*\figheight) -- (-\figlength-\diff,0.75*\figheight) -- (-\figlength-\diff,0) -- (-0.5*\figlength-\diff,-0.5*\figheight) -- cycle;
			\draw (\diff,0) -- (\diff,0.75*\figheight) -- (\figlength/2+\diff,1.25*\figheight) -- (\figlength+\diff,0.75*\figheight) -- (\figlength+\diff,0) -- (0.5*\figlength+\diff,-0.5*\figheight) -- cycle;
			
			\draw[blue] (1.25*\diff,0.25*\diff) -- (1.25*\diff,0.75*\figheight-0.25*\diff) -- (\figlength/2+\diff,1.25*\figheight-0.35*\diff) -- (\figlength+0.75*\diff,0.75*\figheight-0.25*\diff) -- (\figlength+0.75*\diff,0.25*\diff) -- (0.5*\figlength+\diff,-0.5*\figheight+0.35*\diff) -- cycle;
			\draw[blue] (-1.25*\diff,0.25*\diff) -- (-1.25*\diff,0.75*\figheight-0.25*\diff) -- (-\figlength/2-\diff,1.25*\figheight-0.35*\diff) -- (-\figlength-0.75*\diff,0.75*\figheight-0.25*\diff) -- (-\figlength-0.75*\diff,0.25*\diff) -- (-0.5*\figlength-\diff,-0.5*\figheight+0.35*\diff) -- cycle;
			
			\coordinate (pivot) at (-2*\figlength/5-\diff,2*\figheight/4);
			
			\foreach \x in {\figlength/5, 2*\figlength/5,3*\figlength/5,4*\figlength/5}{
				\node[draw, circle, gray, fill=gray, inner sep=\noderadius] at (\x+\diff,0) {};
				\node[draw, circle, fill=black, inner sep=\noderadius] at (-\x-\diff,0) {};
				
				\draw[->, red] (pivot) to [bend left=10] (\x,0.5*\diff);
			}
			\foreach \x in {2*\figlength/6, 3*\figlength/6,4*\figlength/6}{
				\node[draw, circle, gray, fill=gray, inner sep=\noderadius] at (\x+\diff,\figheight/4) {};
				\node[draw, circle, fill=black, inner sep=\noderadius] at (-\x-\diff,\figheight/4) {};
				\draw[->, red] (pivot) to [bend left=10] (\x,0.5*\diff+\figheight/4);
			}
			\foreach \x in {3*\figlength/5}{
				\node[draw, circle, fill=black, inner sep=\noderadius] at (-\x-\diff,2*\figheight/4) {};
			}
			\foreach \x in {2*\figlength/5, 3*\figlength/5}{
				\node[draw, circle, gray, fill=gray, inner sep=\noderadius] at (\x+\diff,2*\figheight/4) {};
				\draw[->, red] (pivot) to [bend left=10] (\x,0.5*\diff+2*\figheight/4);
			}
			\node[draw, circle, gray, fill=gray, inner sep=\noderadius] at (\figlength/2+\diff,3*\figheight/4) {};
			\draw[->, red] (pivot) to [bend left=10] (\figlength/2,0.5*\diff+3*\figheight/4);
			\node[draw, circle, fill=black, inner sep=\noderadius] at (-\figlength/2-\diff,3*\figheight/4) {};

			\node[draw, circle, black, fill=black, inner sep=\noderadius] at (pivot) {};
			\node[draw, circle, cyan, inner sep=2*\noderadius, thick] at (pivot) {};

			\node at (current bounding box.north west) {DG};
		\end{tikzpicture}
        \caption{Couplings of a DG degree of freedom \textcircle{cyan} with all dofs from one of the neighbouring elements.}
		\label{fig:DGcouplings}
	\end{minipage}
    \hfill
	\begin{minipage}[t]{.31\textwidth}
        \centering
		\begin{tikzpicture}
			\draw (-\diff,0) -- (-\diff,0.75*\figheight) -- (-\figlength/2-\diff,1.25*\figheight) -- (-\figlength-\diff,0.75*\figheight) -- (-\figlength-\diff,0) -- (-0.5*\figlength-\diff,-0.5*\figheight) -- cycle;
			\draw (\diff,0) -- (\diff,0.75*\figheight) -- (\figlength/2+\diff,1.25*\figheight) -- (\figlength+\diff,0.75*\figheight) -- (\figlength+\diff,0) -- (0.5*\figlength+\diff,-0.5*\figheight) -- cycle;
			
			\draw[blue] (1.25*\diff,0.25*\diff) -- (1.25*\diff,0.75*\figheight-0.25*\diff) -- (\figlength/2+\diff,1.25*\figheight-0.35*\diff) -- (\figlength+0.75*\diff,0.75*\figheight-0.25*\diff) -- (\figlength+0.75*\diff,0.25*\diff) -- (0.5*\figlength+\diff,-0.5*\figheight+0.35*\diff) -- cycle;
			\draw[blue] (-1.25*\diff,0.25*\diff) -- (-1.25*\diff,0.75*\figheight-0.25*\diff) -- (-\figlength/2-\diff,1.25*\figheight-0.35*\diff) -- (-\figlength-0.75*\diff,0.75*\figheight-0.25*\diff) -- (-\figlength-0.75*\diff,0.25*\diff) -- (-0.5*\figlength-\diff,-0.5*\figheight+0.35*\diff) -- cycle;
			
			\coordinate (pivot) at (-2*\figlength/5-\diff,0);
			
			\foreach \x in {\figlength/5, 2*\figlength/5,3*\figlength/5,4*\figlength/5}{
				\node[draw, circle, gray, fill=gray, inner sep=\noderadius] at (\x+\diff,0) {};
				\node[draw, circle, fill=black, inner sep=\noderadius] at (-\x-\diff,0) {};
				
				\draw[->, red] (pivot) to [bend left=10] (\x,0.5*\diff);
			}
			\foreach \x in {2*\figlength/6, 3*\figlength/6,4*\figlength/6}{
				\node[draw, circle, orange, fill=orange, inner sep=\noderadius] at (\x+\diff,\figheight/4) {};
				\node[draw, circle, orange, fill=orange, inner sep=\noderadius] at (-\x-\diff,\figheight/4) {};
				\draw[->, red] (pivot) to [bend left=10] (\x,0.5*\diff+\figheight/4);
			}
			\foreach \x in {3*\figlength/5, 2*\figlength/5}{
				\node[draw, circle, fill=white, inner sep=\noderadius] at (-\x-\diff,2*\figheight/4) {};
				\node[draw, circle, fill=white, inner sep=\noderadius] at (\x+\diff,2*\figheight/4) {};
			}
			\node[draw, circle, fill=white, inner sep=\noderadius] at (\figlength/2+\diff,3*\figheight/4) {};
			\node[draw, circle, fill=white, inner sep=\noderadius] at (-\figlength/2-\diff,3*\figheight/4) {};

			\node[draw, circle, fill=black, inner sep=\noderadius] at (pivot) {};
			\node[draw, circle, cyan, inner sep=2*\noderadius, thick] at (pivot) {};

			\node at (current bounding box.north west) {TDG};
		\end{tikzpicture}
		\caption{
        A Trefftz DG dof \textcircle{cyan} couples like standard DG dof, however, the number of dofs is reduced by \textcircle{black} for TDG2 and further by all dofs marked \textcircle{orange, fill=orange} for TDG1.}
		\label{fig:TDGcouplings}
	\end{minipage}
	\hfill
	\begin{minipage}[t]{.31\textwidth}
        \centering
		\begin{tikzpicture}
			\draw (-\diff,0) -- (-\diff,0.75*\figheight) -- (-\figlength/2-\diff,1.25*\figheight) -- (-\figlength-\diff,0.75*\figheight) -- (-\figlength-\diff,0) -- (-0.5*\figlength-\diff,-0.5*\figheight) -- cycle;
			\draw (\diff,0) -- (\diff,0.75*\figheight) -- (\figlength/2+\diff,1.25*\figheight) -- (\figlength+\diff,0.75*\figheight) -- (\figlength+\diff,0) -- (0.5*\figlength+\diff,-0.5*\figheight) -- cycle;

			\draw[blue] (0,0) -- (0,0.75*\figheight);
			\draw[blue] (2*\diff+\figlength,0) -- (2*\diff+\figlength,0.75*\figheight);
			\draw[blue] (-2*\diff-\figlength,0) -- (-2*\diff-\figlength,0.75*\figheight);
			
			\draw[blue] (\diff,0.75*\figheight+\diff) -- (0.5*\diff+\figlength/2,1.25*\figheight+0.75*\diff);
			\draw[blue] (-\diff,0.75*\figheight+\diff) -- (-0.5*\diff-\figlength/2,1.25*\figheight+0.75*\diff);
			
			\draw[blue] (\diff,-\diff) -- (0.5*\diff+\figlength/2,-0.5*\figheight-0.75*\diff);
			\draw[blue] (-\diff,-\diff) -- (-0.5*\diff-\figlength/2,-0.5*\figheight-0.75*\diff);
			
			\draw[blue] (1.5*\diff+0.5*\figlength,0.75*\diff+1.25*\figheight) -- (\diff+\figlength,0.75*\figheight+\diff);
			\draw[blue] (1.5*\diff+0.5*\figlength,-0.75*\diff-0.5*\figheight) -- (\diff+\figlength,-\diff);
			\draw[blue] (-1.5*\diff-0.5*\figlength,0.75*\diff+1.25*\figheight) -- (-\diff-\figlength,0.75*\figheight+\diff);
			\draw[blue] (-1.5*\diff-0.5*\figlength,-0.75*\diff-0.5*\figheight) -- (-\diff-\figlength,-\diff);

			\foreach \x in {\figlength/5, 2*\figlength/5,3*\figlength/5,4*\figlength/5}{
				\node[draw, circle, fill=white!90!black, inner sep=\noderadius] at (\x+\diff,0) {};
				\node[draw, circle, fill=white!90!black, inner sep=\noderadius] at (-\x-\diff,0) {};
				
			}
			\foreach \x in {2*\figlength/6, 3*\figlength/6,4*\figlength/6}{
				\node[draw, circle, fill=white!90!black, inner sep=\noderadius] at (\x+\diff,\figheight/4) {};
				\node[draw, circle, fill=white!90!black, inner sep=\noderadius] at (-\x-\diff,\figheight/4) {};
			}
			\foreach \x in {3*\figlength/5}{
				\node[draw, circle, fill=white!90!black, inner sep=\noderadius] at (-\x-\diff,2*\figheight/4) {};
			}
			\foreach \x in {2*\figlength/5, 3*\figlength/5}{
				\node[draw, circle, fill=white!90!black, inner sep=\noderadius] at (\x+\diff,2*\figheight/4) {};
			}
			\node[draw, circle, fill=white!90!black, inner sep=\noderadius] at (\figlength/2+\diff,3*\figheight/4) {};
			\node[draw, circle, fill=white!90!black, inner sep=\noderadius] at (-\figlength/2-\diff,3*\figheight/4) {};
			\node[draw, circle, fill=white!90!black, inner sep=\noderadius] at (-2*\figlength/5-\diff,2*\figheight/4) {};

			\coordinate (pivot) at (0,3*\figheight/5*0.75);
			
			\foreach \x in {-2*\diff-\figlength, 2*\diff+\figlength}{
				\foreach \y in {0.75*\figheight/5, 0.75*2*\figheight/5, 0.75*3*\figheight/5, 0.75*4*\figheight/5}{
					\node[draw, circle, fill=gray, inner sep=0.75*\noderadius] at (\x,\y) {};
				}
			}
			\foreach \y in {0.75*\figheight/5, 0.75*2*\figheight/5, 0.75*4*\figheight/5}{	
				\node[draw, circle, fill=black, inner sep=0.75*\noderadius] at (0,\y) {};
			}
			
			\foreach \idx in {1,2,3,4}{
				\node[draw, circle, fill=gray, inner sep=0.75*\noderadius] at (\diff - \idx*\diff/5 + 0.5*\idx*\diff/5+\idx*\figlength/10, 0.75*\figheight+\diff - 0.75*\idx*\figheight/5 - \idx*\diff/5 + 1.25*\idx*\figheight/5+0.75*\idx*\diff/5) {};
				\node[draw, circle, fill=gray, inner sep=0.75*\noderadius] at (-\diff + \idx*\diff/5 - 0.5*\idx*\diff/5-\idx*\figlength/10, 0.75*\figheight+\diff - 0.75*\idx*\figheight/5 - \idx*\diff/5 + 1.25*\idx*\figheight/5+0.75*\idx*\diff/5) {};
				\node[draw, circle, fill=gray, inner sep=0.75*\noderadius] at (\diff - \idx*\diff/5 + 0.5*\idx*\diff/5+\idx*\figlength/10, -\diff + \idx*\diff/5 - 0.5*\idx*\figheight/5-0.75*\idx*\diff/5) {};
				\node[draw, circle, fill=gray, inner sep=0.75*\noderadius] at (-\diff + \idx*\diff/5 - 0.5*\idx*\diff/5-\idx*\figlength/10, -\diff + \idx*\diff/5 - 0.5*\idx*\figheight/5-0.75*\idx*\diff/5) {};
				
				\node[draw, circle, fill=gray, inner sep=0.75*\noderadius] at (1.5*\diff - \idx*\diff/10 + 0.5*\figlength + \idx*\figlength/10, 0.75*\diff+\idx*\diff/20 + 1.25*\figheight - \idx*\figheight/10) {};
				\node[draw, circle, fill=gray, inner sep=0.75*\noderadius] at (1.5*\diff - \idx*\diff/10 + 0.5*\figlength + \idx*\figlength/10, -0.75*\diff-\idx*\diff/20 - 0.5*\figheight + \idx*\figheight/10) {};
				\node[draw, circle, fill=gray, inner sep=0.75*\noderadius] at (-1.5*\diff + \idx*\diff/10 - 0.5*\figlength - \idx*\figlength/10, -0.75*\diff-\idx*\diff/20 - 0.5*\figheight + \idx*\figheight/10) {};
				\node[draw, circle, fill=gray, inner sep=0.75*\noderadius] at (-1.5*\diff + \idx*\diff/10 - 0.5*\figlength - \idx*\figlength/10, 0.75*\diff+\idx*\diff/20 + 1.25*\figheight - \idx*\figheight/10) {};
			}

			\foreach \x in {-2*\diff-\figlength, 2*\diff+\figlength}{
				\draw[red] (\x,0.75*0.5*\figheight) circle [x radius=\figlength/25, y radius=0.75*0.6*\figheight, rotate=0];
			}
			\foreach \x in {-1.25*\diff-\figlength*0.75, 0.75*\diff+0.25*\figlength}{
				\draw[red] (\x,\figheight+0.875*\diff) circle [x radius=\figlength/25, y radius=0.75*0.6*\figheight, rotate=-56];
				\draw[red] (\x,-3*\diff) circle [x radius=\figlength/25, y radius=0.75*0.6*\figheight, rotate=56];
				\draw[->, red] (pivot) to [bend left=10] (\x, -2.375*\diff);
				\draw[->, red] (pivot) to [bend left=10] (\x, \figheight);
			}
			\foreach \x in {-1.25*\diff-\figlength*0.75, 0.75*\diff+0.25*\figlength}{
				\draw[red] (\x+0.5*\figlength+0.5*\diff,\figheight+0.875*\diff) circle [x radius=\figlength/25, y radius=0.75*0.6*\figheight, rotate=56];
				\draw[red] (\x+0.5*\figlength+0.5*\diff,-3*\diff) circle [x radius=\figlength/25, y radius=0.75*0.6*\figheight, rotate=-56];
				\draw[->, red] (pivot) to [bend right=10] (\x+0.5*\figlength, \figheight);
				\draw[->, red] (pivot) to [bend right=10] (\x+0.5*\figlength, -2.375*\diff);
			}
			\draw[->, red] (pivot) to [bend left=10] (-1.35*\diff-\figlength, 0.75*0.5*\figheight);
			\draw[->, red] (pivot) to [bend right=10] (1.35*\diff+\figlength, 0.75*0.5*\figheight);
			
			\node[draw, circle, fill=blue, inner sep=0.75*\noderadius] at (pivot) {};
			\node[draw, circle, cyan, inner sep=2*\noderadius, thick] at (pivot) {};

			\node at (current bounding box.north west) {HDG};
		\end{tikzpicture}
		\caption{HDG facet dof \textcircle{cyan} couples with all dofs on facets directly adjacent to neighbouring elements.
        }
		\label{fig:HDGcouplings}
	\end{minipage}
\end{figure}

\subsection{DG formulation of the model problem}
We start with the Discontinuous Galerkin (DG) formulation of the model problem. On the mesh $\Th$ we introduce the finite element space of piecewise polynomials of degree $k$:
\begin{subequations}
\begin{equation}
  V_h^k := \{ v \in L^2(\Omega) \mid v|_T \in \mathcal{P}^k(T), \forall T \in \Th \}.
\end{equation}
The symmetric interior penalty DG formulation of the model problem reads, cf. e.g. \cite{di2011mathematical}: Find $u_h \in \Vhk$ such that
\begin{align}
  a_h^\text{DG}(u_h,v_h) &= f_h^\text{DG}(v_h) \quad \forall v \in \Vhk, \text{ with}\\
  a_h^\text{DG}(u,v) &\coloneqq (\nabla u,\nabla v)_{\Th} 
  - (\avg{\nabla u \cdot n}, \jump{v})_{\mathcal{F}_h^{\text{int}}} 
  - (\avg{\nabla v \cdot n}, \jump{u})_{\mathcal{F}_h^{\text{int}}} 
  + (\frac{\lambda}{h} \jump{u}, \jump{v})_{\mathcal{F}_h^{\text{int}}} \\[-1ex]
  &\hspace*{2.23cm}%
   - ({\nabla u \cdot n}, {v})_{\mathcal{F}_h^{\text{bnd}}} \nonumber
 - ({\nabla v \cdot n}, {u})_{\mathcal{F}_h^{\text{bnd}}} 
  + (\frac{\lambda}{h} u,v)_{\mathcal{F}_h^{\text{bnd}}}, \\[-1ex]
  f_h^\text{DG}(v) &\coloneqq 
  - ({\nabla v \cdot n}, u_D)_{\mathcal{F}_h^{\text{bnd}}} 
  + (\frac{\lambda}{h} u_D,v)_{\mathcal{F}_h^{\text{bnd}}}.
\end{align}
\end{subequations}
with $\avg{\cdot}$ and $\jump{\cdot}$ the average and jump operators across facets, $n$ the unit normal vector, $h$ the mesh size, and $\lambda$ a penalty parameter. $\mathcal{F}_h^{\text{int}}$ and $\mathcal{F}_h^{\text{bnd}}$ denote the set of interior and boundary facets, respectively.
The formulation is well known to be stable and convergent for $\lambda \geq c_\lambda k^2$ with a $c_\lambda$ only depending on the shape regularity of the mesh. 
A straightforward static condensation\footnote{see below for further explanation on the term "static condensation"} of degrees of freedom is not possible with the DG formulation. In general, all degrees of freedom of one element are coupled with all other degrees of freedom of a neighbor, cf. \cref{fig:DGcouplings}.

\subsection{Hybrid DG (and Hybrid High Order) formulation of the model problem}
\begin{figure}
	\centering
	\def\diff{1.5mm}
	\def\figheight{2.5cm}
	\def\figlength{4.5cm}
	\def\noderadius{\figlength/3cm}
        \centering
        \begin{tikzpicture}
            \draw (\diff,0) -- (\diff + 0.5*\figlength,0.5*\figheight) -- (\diff,\figheight) -- cycle;
            \draw (-\diff,0) -- (-\diff - 0.5*\figlength,0.5*\figheight) -- (-\diff,\figheight) -- cycle;
            \draw[blue] (-1.25*\diff,0.5*\diff) -- (-0.5*\diff - 0.5*\figlength,0.5*\figheight) -- (-1.25*\diff,\figheight-0.5*\diff) -- cycle;
            \draw[blue] (0,0) -- (0,\figheight);
            \foreach \y in {\figheight/5,2*\figheight/5,4*\figheight/5}{
                \node[draw, circle, fill=black, inner sep=\noderadius] at (-\diff-\figlength/12, \y) {};
            }
            \foreach \y in {0.3*\figheight,0.5*\figheight,0.7*\figheight}{
                \node[draw, circle, fill=black, inner sep=\noderadius] at (-\diff-2*\figlength/12, \y) {};
            }
            \foreach \y in {0.4*\figheight,0.6*\figheight}{
                \node[draw, circle, fill=black, inner sep=\noderadius] at (-\diff-3*\figlength/12, \y) {};
            }
            \node[draw, circle, fill=black, inner sep=\noderadius] at (-\diff-4*\figlength/12, 0.5*\figheight) {};
            
            \draw[red] (-0.25*\figlength-0.75*\diff,0.75*\figheight+1.25*\diff) circle [x radius=\figlength/45, y radius=0.75*0.6*\figheight, rotate=-59];
            \draw[red] (-0.25*\figlength-0.75*\diff,0.25*\figheight-1.25*\diff) circle [x radius=\figlength/45, y radius=0.75*0.6*\figheight, rotate=59];
            \draw[red] (0,0.5*\figheight) circle [x radius=\figlength/45, y radius=0.75*0.6*\figheight, rotate=0];

            \foreach \y in {\figheight/5,2*\figheight/5,3*\figheight/5,4*\figheight/5}{
                \node[draw, circle, fill=black, inner sep=\noderadius] at (\diff+\figlength/12, \y) {};
            }
            \foreach \y in {0.3*\figheight,0.5*\figheight,0.7*\figheight}{
                \node[draw, circle, fill=black, inner sep=\noderadius] at (\diff+2*\figlength/12, \y) {};
            }
            \foreach \y in {0.4*\figheight,0.6*\figheight}{
                \node[draw, circle, fill=black, inner sep=\noderadius] at (\diff+3*\figlength/12, \y) {};
            }
            \node[draw, circle, fill=black, inner sep=\noderadius] at (\diff+4*\figlength/12, 0.5*\figheight) {};
            
            \coordinate (pivot) at (-\diff-\figlength/12,0.6*\figheight);

            \foreach \y in {0.2*\figheight,0.4*\figheight,0.6*\figheight,0.8*\figheight}{
                \node[draw, circle, fill=gray, inner sep=\noderadius] at (0,\y) {};
            }
            \draw[->, red] (pivot) to [bend right=10] (-0.25*\figlength,0.25*\figheight);
            \draw[->, red] (pivot) to [bend left=10] (-0.25*\figlength,0.75*\figheight);
            \draw[->, red] (pivot) to [bend right=10] (-1.25*\diff,0.5*\figheight);
            
            \draw[blue] (-\diff-0.5*\figlength,0.5*\figheight+\diff) -- (-\diff,\figheight+\diff);
            \draw[blue] (-\diff-0.5*\figlength,0.5*\figheight-\diff) -- (-\diff,-\diff);
            \draw[white] (\diff+0.5*\figlength,0.5*\figheight-\diff) -- (\diff,-\diff);
            \draw[white] (\diff+0.5*\figlength,0.5*\figheight+\diff) -- (\diff,\figheight+\diff);
            
            \foreach \idx in {1,2,3,4}{
                \node[draw, circle, fill=gray, inner sep=\noderadius] at (-\diff-\idx*\figlength/10,\idx*\figheight/10-\diff) {};
                \node[draw, circle, fill=gray, inner sep=\noderadius] at (-\diff-\idx*\figlength/10,\figheight-\idx*\figheight/10+\diff) {};
            }
            
            \node[draw, circle, fill=black, inner sep=\noderadius] at (pivot) {};
            \node[draw, circle, cyan, inner sep=2*\noderadius, thick] at (pivot) {};
            \node[below]
            at (current bounding box.north west) {HDG \textcolor{gray}{(wo)}};
        \end{tikzpicture}
        \hspace{5em}
        \begin{tikzpicture}
            \draw (\diff,0) -- (\diff + 0.5*\figlength,0.5*\figheight) -- (\diff,\figheight) -- cycle;
            \draw (-\diff,0) -- (-\diff - 0.5*\figlength,0.5*\figheight) -- (-\diff,\figheight) -- cycle;
            \draw[blue] (0,0) -- (0,\figheight);
            \draw[blue] (-\diff-0.5*\figlength,0.5*\figheight+\diff) -- (-\diff,\figheight+\diff);
            \draw[blue] (-\diff-0.5*\figlength,0.5*\figheight-\diff) -- (-\diff,-\diff);
            \draw[blue] (\diff+0.5*\figlength,0.5*\figheight-\diff) -- (\diff,-\diff);
            \draw[blue] (\diff+0.5*\figlength,0.5*\figheight+\diff) -- (\diff,\figheight+\diff);
            \foreach \y in {\figheight/5,2*\figheight/5,4*\figheight/5}{
                \node[draw, circle, fill=white!90!black, inner sep=\noderadius] at (-\diff-\figlength/12, \y) {};
            }
            \foreach \y in {0.3*\figheight,0.5*\figheight,0.7*\figheight}{
                \node[draw, circle, fill=white!90!black, inner sep=\noderadius] at (-\diff-2*\figlength/12, \y) {};
            }
            \foreach \y in {0.4*\figheight,0.6*\figheight}{
                \node[draw, circle, fill=white!90!black, inner sep=\noderadius] at (-\diff-3*\figlength/12, \y) {};
            }
            \node[draw, circle, fill=white!90!black, inner sep=\noderadius] at (-\diff-4*\figlength/12, 0.5*\figheight) {};
            \node[draw, circle, fill=white!90!black, inner sep=\noderadius] at (-\diff-\figlength/12, 3*\figheight/5) {};

            \foreach \y in {\figheight/5,2*\figheight/5,3*\figheight/5,4*\figheight/5}{
                \node[draw, circle, fill=white!90!black, inner sep=\noderadius] at (\diff+\figlength/12, \y) {};
            }
            \foreach \y in {0.3*\figheight,0.5*\figheight,0.7*\figheight}{
                \node[draw, circle, fill=white!90!black, inner sep=\noderadius] at (\diff+2*\figlength/12, \y) {};
            }
            \foreach \y in {0.4*\figheight,0.6*\figheight}{
                \node[draw, circle, fill=white!90!black, inner sep=\noderadius] at (\diff+3*\figlength/12, \y) {};
            }
            \node[draw, circle, fill=white!90!black, inner sep=\noderadius] at (\diff+4*\figlength/12, 0.5*\figheight) {};
            
            \coordinate (pivot) at (0, 3*\figheight/5);

            \foreach \y in {0.2*\figheight,0.4*\figheight,0.8*\figheight}{
                \node[draw, circle, fill=black, inner sep=\noderadius] at (0,\y) {};
            }
            \foreach \idx in {1,2,3,4}{
                \node[draw, circle, fill=gray, inner sep=\noderadius] at (\diff+\idx*\figlength/10,\idx*\figheight/10-\diff) {};
                \node[draw, circle, fill=gray, inner sep=\noderadius] at (-\diff-\idx*\figlength/10,\idx*\figheight/10-\diff) {};
                \node[draw, circle, fill=gray, inner sep=\noderadius] at (\diff+\idx*\figlength/10,\figheight-\idx*\figheight/10+\diff) {};
                \node[draw, circle, fill=gray, inner sep=\noderadius] at (-\diff-\idx*\figlength/10,\figheight-\idx*\figheight/10+\diff) {};
                
            }
        
            \draw[->, red] (pivot) to [bend right=10] (-0.25*\figlength,0.25*\figheight);
            \draw[->, red] (pivot) to [bend left=10] (-0.25*\figlength,0.75*\figheight);
            \draw[->, red] (pivot) to [bend left=10] (0.25*\figlength,0.25*\figheight);
            \draw[->, red] (pivot) to [bend right=10] (0.25*\figlength,0.75*\figheight);

            \draw[red] (-0.25*\figlength-0.75*\diff,0.75*\figheight+1.25*\diff) circle [x radius=\figlength/45, y 	radius=0.75*0.6*\figheight, rotate=-59];
            \draw[red] (-0.25*\figlength-0.75*\diff,0.25*\figheight-1.25*\diff) circle [x radius=\figlength/45, y 	radius=0.75*0.6*\figheight, rotate=59];
            \draw[red] (0.25*\figlength+0.75*\diff,0.75*\figheight+1.25*\diff) circle [x radius=\figlength/45, y 	radius=0.75*0.6*\figheight, rotate=59];
            \draw[red] (0.25*\figlength+0.75*\diff,0.25*\figheight-1.25*\diff) circle [x radius=\figlength/45, y 	radius=0.75*0.6*\figheight, rotate=-59];
            
            \node[draw, circle, fill=black, inner sep=\noderadius] at (pivot) {};
            \node[draw, circle, cyan, inner sep=2*\noderadius, thick] at (pivot) {};
            \node[below]
            at (current bounding box.north west) {HDG \textcolor{gray}{(w)}};
        \end{tikzpicture} 

    \caption{Comparison of Hybrid DG dof of order $3$ with (w) and without (wo) static condensation. An element dof \textcircle{cyan} couples with its own elements dof and those of neighbouring facets, but different volume elements no longer interact directly in the resulting matrix. Further, by applying static condensation element dofs can be removed from the global system, remaining only with facet dofs and couplings to their adjacent counterparts.}
    \label{fig:staticcond}
\end{figure}
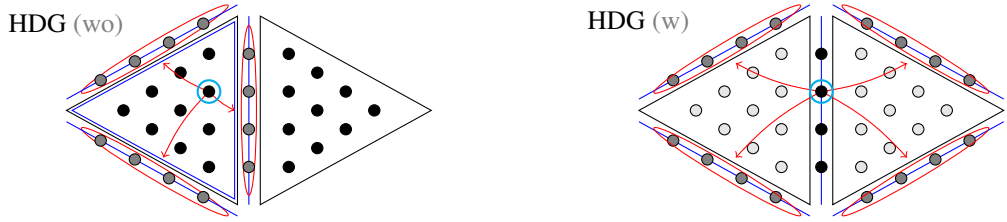

The idea of Hybrid DG methods is to introduce additional facet unknowns that allow the removal of the direct couplings between element unknowns of neighboring elements. To this end, we introduce the space of facet unknowns and the product space of element and facet unknowns, both with homogeneous and prescribed Dirichlet boundary data:
\begin{subequations}
\begin{align}
  \Fhkz &:= \{ v_F \in L^2(\mathcal{F}_h) \mid v_F|_F \in \mathcal{P}^{k}(F) ~\forall F \in \mathcal{F}_h^{\text{int}}, v_F|_F = 0 ~\forall F \in \mathcal{F}_h^{\text{bnd}} \},\!\!\!\!  &\Whkz &:= \Vhk \times \Fhkz,
  \\
  \Fhkd &:= \{ v_F \in L^2(\mathcal{F}_h) \mid v_F|_F \in \mathcal{P}^{k}(F) ~\forall F \in \mathcal{F}_h^{\text{int}}, v_F|_F = \Pi_F u_D ~\forall F \in \mathcal{F}_h^{\text{bnd}} \},\!\!\!\! &\Whkd &:= \Vhk \times \Fhkd,
\end{align}
where $\Pi_F$ denotes the $L^2$-projection onto the facet polynomials $\mathcal{P}^k(F)$.
A corresponding hybrid DG formulation of the model problem reads: 
Find $u_h=(u_T,u_F) \in \Whkd$ such that
  \begin{align}
    a_h^\text{HDG}(u_h,v_h) &= 0 \quad \forall v_h = (v_T,v_F) \in \Whkz, \text{ with}\\
    a_h^\text{HDG}(u_h,v_h) &\coloneqq (\nabla u_T,\nabla v_T)_{\Th} 
    - (\nabla u_T \cdot n, \jump{v_h})_{\partial\mathcal{T}_h} 
    - (\nabla v_T \cdot n, \jump{u_h})_{\partial\mathcal{T}_h} 
    + (\frac{\lambda}{h} \jump{u_h}, \jump{v_h})_{\partial\mathcal{T}_h}
  \end{align}
  \end{subequations}
where $\partial\mathcal{T}_h$ denotes the set of all element boundaries and $\lambda$ is chosen as for the DG formulation. 
\paragraph{Static condensation}
Now, it is important to note that the element unknowns $u_T$ can be determined locally solely based on the facet variable $u_F$ on the adjacent facets. This allows elimination of the element volume unknowns from the system of equations by a Schur complement strategy known as \emph{static condensation}. The resulting system of equations is then formulated only in terms of the facet unknowns $u_F$ which are -- especially for higher order methods -- significantly less than the element unknowns. 
Instead of number of element unknowns that scale with $k^d$ the facet unknowns per facet only scale with $k^{d-1}$. In \cref{fig:staticcond} a sketch of the coupling pattern before and after static condensation for the HDG method on triangles is shown. 

The coupling pattern for the HDG method (after static condensation) on hexagons in comparison to the DG method is shown in \cref{fig:HDGcouplings}.

\paragraph{Hybrid High Order (HHO) methods}
An improvement in the efficiency of the HDG method for diffusion-dominated problems can be achieved if facet unknowns are chosen in $F_h^{k-1}$ instead of $F_h^k$. To achieve this in a consistent manner the discrete variational form is adjusted or post-processing techniques have to be applied.

Systematically, for polytopes, the framework of Hybrid High Order (HHO) methods achieves this. Instead of going into detail here, we refer to \cite{di2015hybrid,di2016review} for a detailed discussion of the HHO method. The HHO method is however closely related to the HDG method which can also be brought into the same structure in terms of the facet unknowns, cf. \cite{cockburn2016bridging}. Key to this is either a modification of the jump stabilization with the \emph{projected jumps} approach, cf. \cite{lehrenfeld2016high}, or a post-processing step to lift an order $k-1$ approximation in the interior to an accurate order $k$ approximation. On polytopal meshes, this requires some additional care but is feasible with the concept of M-decompositions, cf. \cite{cockburn2017superconvergence1,cockburn2017superconvergence2,cockburn2018introduction}.
In the remainder of this paper, we will focus on the coupling structure of methods and denote the HDG method with a reduced order on the facet as HHO without considering if and when the HHO method may not be applicable.

\subsection{Trefftz DG formulation of the model problem}
An alternative approach to reduce the number of unknowns in the DG method is the Trefftz DG method. Here, the basis functions are chosen to be solutions of the (homogeneous) PDE on each element. The Trefftz space is defined as
\begin{subequations}
\begin{align} \label{eq:trefftzspace}
  \TT_h^k=\{v\in \Vhk,\ - \Delta v|_T=0, \quad T\in\Th\}\subset \Vhk.
\end{align}
Note that for different PDEs the Trefftz space will be different and possibly the PDE constraint within the space may be relaxed as in the quasi- and weak Trefftz methods \cite{LS_IJMNE_2023,IGMS_MC_2021,LLS_NM_2024,perinati2023quasi}.
With that choice the discrete variational formulation reads:
Find $u_h \in \TT_h^k$ such that
\begin{align}
  a_h^\text{TDG}(u_h,v_h) &= f_h^\text{TDG}(v_h) \quad \forall v \in \TT_h^k
\end{align}
\end{subequations}
where $a_h^\text{TDG}(u_h,v_h)=a_h^\text{DG}(u_h,v_h)$ and $f_h^\text{TDG}(v_h)=f_h^\text{DG}(v_h)$, i.e. we use the same bi- and linear form as in the DG method.

Let us briefly discuss the dimension reduction of the Trefftz DG method. The number of degrees of freedom on one element for the Trefftz DG method is given by the dimension of the kernel of the Laplace operator which is the difference between the dimension of the space of piecewise polynomials of degree $k$, $\mathcal{P}^k(T)$, and the dimension of the space of polynomials of degree $k-2$, $\mathcal{P}^{k-2}(T) = \Delta \mathcal{P}^k(T)$. This leads to a scaling of the unknowns on each element with $k^{d-1}$ instead of $k^d$ as in the DG method.
In contrast to the HDG method, the scaling reduction is achieved while keeping the unknowns on the element instead of moving them to the facets, cf. \cref{fig:TDGcouplings}.

Note that for different differential operators the kernel (and its dimension) might be different. Especially, for first order operators, as for instance advection $\mathbf{w}\cdot \nabla$ (where $\mathbf{w}$ is a suitable flow field), the dimension of the kernel is even smaller, it is the difference in dimension of the space of piecewise polynomials of degree $k$ and degree $k-1$. In the following, we will abbreviate the Trefftz DG method as $\text{TDG1}$ for first order and $\text{TDG2}$ for second order (scalar) operators.

The Trefftz DG method has traditionally been used only for a few specific problems, such as the Helmholtz or the Laplace equation, where the Trefftz space is the space of plane waves. However, with the introduction of Quasi-Trefftz \cite{imbert2014generalized, IGMS_MC_2021, perinati2023quasi} and embedded Trefftz methods \cite{LS_IJMNE_2023}, the Trefftz DG method can be applied to a wider range of problems, including problems with variable coefficients and inhomogeneous r.h.s.. Although the theory of these methods is still under development, we assume in the following that the Trefftz DG method can be applied and achieves the optimal dimension reduction as described above.

\subsection{Virtual Element Method formulation of the model problem}
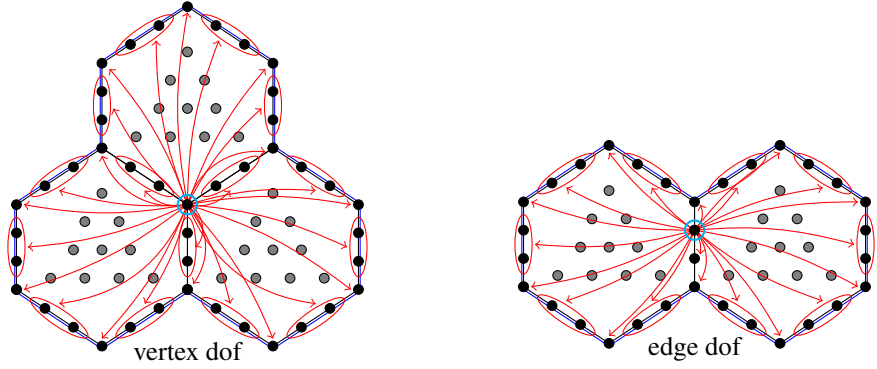
\begin{figure}
	\centering
	\def\figheight{1.5cm}
	\def\figlength{1.5*\figheight}
	\def\noderadius{\figlength/1.75cm}
	\def\diff{1.5mm}
	\begin{tikzpicture}
		\draw (0,0) -- (0,0.75*\figheight) -- (-\figlength/2,1.25*\figheight) -- (-\figlength,0.75*\figheight) -- (-\figlength,0) -- (-0.5*\figlength,-0.5*\figheight) -- cycle;
		\draw (0,0) -- (0,0.75*\figheight) -- (\figlength/2,1.25*\figheight) -- (\figlength,0.75*\figheight) -- (\figlength,0) -- (0.5*\figlength,-0.5*\figheight) -- cycle;
		\draw (0,0.75*\figheight) -- (\figlength/2,1.25*\figheight) -- (\figlength/2,2*\figheight) -- (0,2.5*\figheight) -- (-0.5*\figlength,2*\figheight) -- (-0.5*\figlength,1.25*\figheight) -- cycle;
		
		\draw[blue] (-0.15*\diff-\figlength,-0.15*\diff) -- (-0.15*\diff-\figlength,0.75*\figheight+0.125*\diff) -- (-0.15*\diff-\figlength/2, 1.25*\figheight+0.125*\diff)-- (-0.15*\diff-\figlength/2, 2*\figheight+0.125*\diff)-- (0, 2.5*\figheight+0.2*\diff) -- (0.15*\diff+\figlength/2, 2*\figheight+0.125*\diff)-- (0.15*\diff+\figlength/2, 1.25*\figheight+0.125*\diff) -- (0.15*\diff+\figlength, 0.75*\figheight+0.15*\diff) -- (0.15*\diff+\figlength, -0.15*\diff) -- (\figlength/2, -0.5*\figheight-0.25*\diff) -- (0,-0.25*\diff) -- (-\figlength/2, -0.5*\figheight-0.25*\diff) -- cycle;

		\coordinate (pivot) at (0,0.75*\figheight) {};
		
		\foreach \x in {\figlength/5, 2*\figlength/5,3*\figlength/5,4*\figlength/5}{
			\node[draw, circle, fill=gray, inner sep=\noderadius] at (\x,\diff) {};
			\node[draw, circle, fill=gray, inner sep=\noderadius] at (-\x,\diff) {};
			\node[draw, circle, fill=gray, inner sep=\noderadius] at (\x-\figlength/2,\diff+1.25*\figheight) {};

		}
		\foreach \x in {2*\figlength/6, 3*\figlength/6,4*\figlength/6}{
			\node[draw, circle, fill=gray, inner sep=\noderadius] at (\x,\figheight/4+\diff) {};
			\node[draw, circle, fill=gray, inner sep=\noderadius] at (\x-\figlength/2,\figheight/4+\diff+1.25*\figheight) {};
			\node[draw, circle, fill=gray, inner sep=\noderadius] at (-\x,\figheight/4+\diff) {};
		}
		\foreach \x in {3*\figlength/5, 2*\figlength/5}{
			\node[draw, circle, fill=gray, inner sep=\noderadius] at (\x,2*\figheight/4+\diff) {};
			\node[draw, circle, fill=gray, inner sep=\noderadius] at (\x-\figlength/2,2*\figheight/4+\diff+1.25*\figheight) {};
			\node[draw, circle, fill=gray, inner sep=\noderadius] at (-\x,2*\figheight/4+\diff) {};
		}
		\node[draw, circle, fill=gray, inner sep=\noderadius] at (\figlength/2,3*\figheight/4+\diff) {};
		\node[draw, circle, fill=gray, inner sep=\noderadius] at (0,3*\figheight/4+\diff+1.25*\figheight) {};
		\node[draw, circle, fill=gray, inner sep=\noderadius] at (-\figlength/2,3*\figheight/4+\diff) {};
		
		\foreach \idx in {0,1,2}{
			\node[draw, circle, fill=black, inner sep=\noderadius] at (\figlength,0.75*\idx*\figheight/3) {};
			\foreach \x in {-\figlength, 0}{
				\node[draw, circle, fill=black, inner sep=\noderadius] at (\x,0.75*\idx*\figheight/3) {};
				\node[draw, circle, fill=black, inner sep=\noderadius] at (\x+0.5*\figlength,0.75*\idx*\figheight/3+1.25*\figheight) {};

			}
			\foreach \y in {2*\figheight+\idx*0.5*\figheight/3, 1.25*\figheight-\idx*0.5*\figheight/3}{
				\node[draw, circle, fill=black, inner sep=\noderadius] at (-0.5*\figlength+0.5*\idx*\figlength/3,\y) {};
			}
			\foreach \y in {0,2.5*\figheight}{
				\node[draw, circle, fill=black, inner sep=\noderadius] at (0.5*\idx*\figlength/3,\y-\idx*0.5*\figheight/3) {};
			}
			\foreach \x in {-\figlength+0.5*\idx*\figlength/3, 0.5*\idx*\figlength/3}{
				\node[draw, circle, fill=black, inner sep=\noderadius] at (\x,0.75*\figheight+\idx*0.5*\figheight/3) {};			
			}
			\foreach \x in {-1,1}{
				\node[draw, circle, fill=black, inner sep=\noderadius] at (\x*0.5*\figlength+0.5*\idx*\figlength/3,-0.5*\figheight+\idx*0.5*\figheight/3) {};
			}
			\node[draw, circle, fill=black, inner sep=\noderadius] at (-\figlength+0.5*\idx*\figlength/3,-\idx*0.5*\figheight/3) {};
			\node[draw, circle, fill=black, inner sep=\noderadius] at (\figlength/2+0.5*\idx*\figlength/3,1.25*\figheight-\idx*0.5*\figheight/3) {};

			\draw[red] (-0.75*\figlength+\idx*0.5*\figlength,-0.25*\figheight+\idx*1.25*\figheight) circle [x radius=\figlength/20, y radius=0.75*0.4*\figheight, rotate=56];
			\draw[red] (0.75*\figlength-\idx*0.5*\figlength,-0.25*\figheight+\idx*1.25*\figheight) circle [x radius=\figlength/20, y radius=0.75*0.4*\figheight, rotate=-56];
			\draw[red] (-1*\figlength+\idx*\figlength,0.75*\figheight/2) circle [x radius=\figlength/20, y radius=0.75*0.35*\figheight, rotate=0];
		}
		\node[draw, circle, fill=black, inner sep=\noderadius] at (\figlength,0.75*\figheight) {};
		\node[draw, circle, fill=black, inner sep=\noderadius] at (\figlength/2,2*\figheight) {};

		\foreach \idx in {0,1}{
			\draw[red] (-0.75*\figlength+\idx*\figlength/2,1*\figheight-\idx*1.25*\figheight) circle [x radius=\figlength/20, y radius=0.75*0.4*\figheight, rotate=-56];
			\draw[red] (0.25*\figlength+\idx*\figlength/2,-0.25*\figheight+\idx*1.25*\figheight) circle [x radius=\figlength/20, y radius=0.75*0.4*\figheight, rotate=56];
			\draw[red] (-0.5*\figlength+\idx*\figlength,1.25*\figheight+0.75*\figheight/2) circle [x radius=\figlength/20, y radius=0.75*0.35*\figheight, rotate=0];

			\draw[->,red] (pivot) to [bend left=10] (-0.25*\figlength-\idx*\figlength/2,-0.25*\figheight+\idx*1.25*\figheight+1.25*\diff-2*\idx*1.25*\diff);
			\draw[->,red] (pivot) to [bend left=10] (-0.5*\figlength+\idx*\figlength+\diff-2*\idx*\diff,1.6*\figheight);
			\draw[->,red] (pivot) to [bend left=10] (-0.25*\figlength+\idx*\figlength/2+\diff-2*\idx*\diff,2.25*\figheight-\diff/2);
			\draw[->,red] (pivot) to [bend left=30] (-0.25*\figlength-\idx*\figlength/2,1*\figheight-\idx*1.25*\figheight-1.25*\diff+2*\idx*1.25*\diff);
			\draw[->,red] (pivot) to [bend right=30] (0.25*\figlength+\idx*\figlength/2,1*\figheight-\idx*1.25*\figheight-1.25*\diff+2*\idx*1.25*\diff);
			\draw[->,red] (pivot) to [bend left=10] (0.25*\figlength+\idx*\figlength/2,-0.25*\figheight+\idx*1.25*\figheight+1.25*\diff-2*\idx*1.25*\diff);
			\draw[->,red] (pivot) to [bend left=10] (-1*\figlength+\diff+2*\idx*\figlength-\idx*2*\diff,0.75*\figheight/2);
			\draw[->,red] (pivot) to [bend left=10] (-1*\figlength+2*\diff/3+2*\idx*\figlength-\idx*2*2*\diff/3,0.75*\figheight);
			\draw[->,red] (pivot) to [bend left=10] (-1*\figlength+2*\diff/3+2*\idx*\figlength-\idx*2*2*\diff/3,0);
			\draw[->,red] (pivot) to [bend left=10] (-0.5*\figlength+\idx*\figlength,-0.5*\figheight+2*\diff/3);
			\draw[->,red] (pivot) to [bend left=40] (-0.5*\figlength+\idx*\figlength-\idx*2*\diff/3,1.25*\figheight-2*\diff/3+\idx*\diff/3);
			\draw[->,red] (pivot) to [bend left=10] (\idx*\figlength/2-\idx*2*\diff/3,2.5*\figheight-\diff/2-\idx*\figheight/2);
		}
		\draw[->,red] (pivot) to [bend right=10] (-\figlength/2+2*\diff/3,2.5*\figheight-\diff/2-\figheight/2);
		\draw[->,red] (pivot) to [bend left=60] (2*\diff/3,0.75*\figheight-0.75*\figheight/2);
		\draw[->,red] (pivot) to [bend left=40] (0.5*\diff,0.5*\diff);

		\node[draw, circle, fill=black, inner sep=\noderadius] at (pivot) {};
        \node[draw, circle, cyan, inner sep=2*\noderadius, thick] at (pivot) {};

		\node at (current bounding box.south) {vertex dof};
	\end{tikzpicture}
    \hspace{5em}
	\begin{tikzpicture}
		\draw (0,0) -- (0,0.75*\figheight) -- (-\figlength/2,1.25*\figheight) -- (-\figlength,0.75*\figheight) -- (-\figlength,0) -- (-0.5*\figlength,-0.5*\figheight) -- cycle;
		\draw (0,0) -- (0,0.75*\figheight) -- (\figlength/2,1.25*\figheight) -- (\figlength,0.75*\figheight) -- (\figlength,0) -- (0.5*\figlength,-0.5*\figheight) -- cycle;
		
		\draw[blue] (-0.15*\diff-\figlength,-0.15*\diff) -- (-0.15*\diff-\figlength,0.75*\figheight+0.125*\diff) -- (-\figlength/2, 1.25*\figheight+0.25*\diff)-- (0, 0.75*\figheight+0.25*\diff) -- (\figlength/2, 1.25*\figheight+0.25*\diff) -- (0.15*\diff+\figlength, 0.75*\figheight+0.125*\diff) -- (0.15*\diff+\figlength, -0.15*\diff) -- (\figlength/2, -0.5*\figheight-0.25*\diff) -- (0,-0.25*\diff) -- (-\figlength/2, -0.5*\figheight-0.25*\diff) -- cycle;

		\coordinate (pivot) at (0,0.75*2*\figheight/3) {};
		
		\foreach \x in {\figlength/5, 2*\figlength/5,3*\figlength/5,4*\figlength/5}{
			\node[draw, circle, fill=gray, inner sep=\noderadius] at (\x,\diff) {};
			\node[draw, circle, fill=gray, inner sep=\noderadius] at (-\x,\diff) {};

		}
		\foreach \x in {2*\figlength/6, 3*\figlength/6,4*\figlength/6}{
			\node[draw, circle, fill=gray, inner sep=\noderadius] at (\x,\figheight/4+\diff) {};
			\node[draw, circle, fill=gray, inner sep=\noderadius] at (-\x,\figheight/4+\diff) {};
		}
		\foreach \x in {3*\figlength/5, 2*\figlength/5}{
			\node[draw, circle, fill=gray, inner sep=\noderadius] at (\x,2*\figheight/4+\diff) {};
			\node[draw, circle, fill=gray, inner sep=\noderadius] at (-\x,2*\figheight/4+\diff) {};
		}
		\node[draw, circle, fill=gray, inner sep=\noderadius] at (\figlength/2,3*\figheight/4+\diff) {};
		\node[draw, circle, fill=gray, inner sep=\noderadius] at (-\figlength/2,3*\figheight/4+\diff) {};

		\foreach \idx in {0,1}{
			\draw[red] (-0.75*\figlength+\idx*\figlength/2,-0.25*\figheight+\idx*1.25*\figheight) circle [x radius=\figlength/20, y radius=0.75*0.4*\figheight, rotate=56];
			\draw[red] (-0.75*\figlength+\idx*\figlength/2,1*\figheight-\idx*1.25*\figheight) circle [x radius=\figlength/20, y radius=0.75*0.4*\figheight, rotate=-56];
			\draw[red] (0.25*\figlength+\idx*\figlength/2,1*\figheight-\idx*1.25*\figheight) circle [x radius=\figlength/20, y radius=0.75*0.4*\figheight, rotate=-56];
			\draw[red] (0.25*\figlength+\idx*\figlength/2,-0.25*\figheight+\idx*1.25*\figheight) circle [x radius=\figlength/20, y radius=0.75*0.4*\figheight, rotate=56];
			\draw[red] (-1*\figlength+2*\idx*\figlength,0.75*\figheight/2) circle [x radius=\figlength/20, y radius=0.75*0.35*\figheight, rotate=0];
			\draw[->,red] (pivot) to [bend left=10] (-0.25*\figlength-\idx*\figlength/2,-0.25*\figheight+\idx*1.25*\figheight+1.25*\diff-2*\idx*1.25*\diff);
			\draw[->,red] (pivot) to [bend left=10] (-0.25*\figlength-\idx*\figlength/2,1*\figheight-\idx*1.25*\figheight-1.25*\diff+2*\idx*1.25*\diff);
			\draw[->,red] (pivot) to [bend right=10] (0.25*\figlength+\idx*\figlength/2,1*\figheight-\idx*1.25*\figheight-1.25*\diff+2*\idx*1.25*\diff);
			\draw[->,red] (pivot) to [bend left=10] (0.25*\figlength+\idx*\figlength/2,-0.25*\figheight+\idx*1.25*\figheight+1.25*\diff-2*\idx*1.25*\diff);
			\draw[->,red] (pivot) to [bend left=10] (-1*\figlength+\diff+2*\idx*\figlength-\idx*2*\diff,0.75*\figheight/2);
			\draw[->,red] (pivot) to [bend left=10] (-1*\figlength+2*\diff/3+2*\idx*\figlength-\idx*2*2*\diff/3,0.75*\figheight);
			\draw[->,red] (pivot) to [bend left=10] (-1*\figlength+2*\diff/3+2*\idx*\figlength-\idx*2*2*\diff/3,0);
			\draw[->,red] (pivot) to [bend left=10] (-0.5*\figlength+\idx*\figlength,-0.5*\figheight+2*\diff/3);
			\draw[->,red] (pivot) to [bend left=10] (-0.5*\figlength+\idx*\figlength,1.25*\figheight-2*\diff/3);
		}
		
		\foreach \idx in {0,1,2}{
			
			\node[draw, circle, fill=black, inner sep=\noderadius] at (\figlength,0.75*\idx*\figheight/3) {};
			\foreach \x in {-\figlength, 0}{
				\node[draw, circle, fill=black, inner sep=\noderadius] at (\x,0.75*\idx*\figheight/3) {};

			}
			\node[draw, circle, fill=black, inner sep=\noderadius] at (0.5*\idx*\figlength/3,-\idx*0.5*\figheight/3) {};
			\foreach \x in {-\figlength+0.5*\idx*\figlength/3, 0.5*\idx*\figlength/3}{
				\node[draw, circle, fill=black, inner sep=\noderadius] at (\x,0.75*\figheight+\idx*0.5*\figheight/3) {};			
			}
			\foreach \x in {-1,1}{
				\node[draw, circle, fill=black, inner sep=\noderadius] at (\x*0.5*\figlength+0.5*\idx*\figlength/3,-0.5*\figheight+\idx*0.5*\figheight/3) {};
			}
			\node[draw, circle, fill=black, inner sep=\noderadius] at (-\figlength/2+0.5*\idx*\figlength/3,1.25*\figheight-\idx*0.5*\figheight/3) {};
			\node[draw, circle, fill=black, inner sep=\noderadius] at (-\figlength+0.5*\idx*\figlength/3,-\idx*0.5*\figheight/3) {};
			\node[draw, circle, fill=black, inner sep=\noderadius] at (\figlength/2+0.5*\idx*\figlength/3,1.25*\figheight-\idx*0.5*\figheight/3) {};
		}
		\node[draw, circle, fill=black, inner sep=\noderadius] at (\figlength,0.75*\figheight) {};
		\draw[->,red] (pivot) to [bend left=35] (0.5*\diff,0.5*\diff);
		\draw[->,red] (pivot) to [bend left=20] (0.5*\diff,0.75*\figheight/4+\diff);
		\draw[->,red] (pivot) to [bend right=35] (0.5*\diff,0.75*\figheight-0.25*\diff);
		
		\node[draw, circle, fill=black, inner sep=\noderadius] at (pivot) {};
        \node[draw, circle, cyan, inner sep=2*\noderadius, thick] at (pivot) {};

		\node at (current bounding box.south) {edge dof};
	\end{tikzpicture}
	\caption{All types of dofs couple with all other dofs that share a common element. Element dofs do not couple with neighboring element dofs and can be condensed out by a Schur complement strategy.} \label{fig:VEMcouplings}
\end{figure}

The Virtual Element Method (VEM), introduced in \cite{vembasicprinciples, vemhitchhiker}, is a generalization of the conforming finite element method to very general polytopal meshes.
The method gets its name from the fact that the conforming basis functions are not explicitly known (unless one uses lightning VEM \cite{trezzizerbinati2023}), but are only defined by their degrees of freedom. The VE space is defined as
\begin{subequations}
\begin{align}
    \VEM^k = \{ v \in H^1(\Omega) : v|_{F}\in\mathcal{P}^k(F),\ \forall F\in\partial T,\ v|_{\partial T}\in C^0(\partial T),\ \Delta v|_T\in\mathcal{P}^{k-2}(T), \quad \forall T \in \Th \},
\end{align}
with $\VEM_0^k$ and $\VEM_D^k$ the subspaces with homogeneous or Dirichlet boundary conditions, respectively.
Note that the last condition, involving the Laplacian, is not to be confused with a Trefftz-type condition, as it is independent of the differential operator of the PDE, but rather required to define the degrees of freedom.
Degrees of freedom for the space are the values at the vertices, $k-1$ values on each edge, and the moments up to order $k-2$ on each volume element and also on faces in 3D, see also \cref{tab:ndof} in \cref{sec:measures} below.
The VEM formulation of the model problem reads: Find $u_h \in \VEM^k_{D}$ such that
\begin{align}
    a_h^\text{VEM}(\Pi_h u_h,\Pi_h v_h) + s_h^\text{VEM}(u_h-\Pi_h u_h, v_h-\Pi_h v_h) = f_h^\text{VEM}(v_h) \quad \forall v_h \in \VEM^k_0,
\end{align}
\end{subequations}
where $a_h^\text{VEM}(\cdot,\cdot)=a(\cdot,\cdot)$, $\Pi_h$ is a projection onto the space of piecewise polynomials of degree $k$ and $s_h^\text{VEM}$ is a stabilization term.
Boundary conditions can be enforced directly using the degrees of freedom on the boundary or by a penalty term.
The point of the projection is to ensure that the bilinear form is computable on the VEM space, whereas the stabilization term is used to ensure the well-posedness of the method.
See e.g. \cite{vembasicprinciples, vembook} for possible choices for the projection and the stabilization term and the theoretical background.
For details on the implementation of the VEM method see e.g. \cite{vemhitchhiker, vembook}.

Note that the couplings in $a_h^{\text{VEM}}(\cdot,\cdot)$ and $s_h^{\text{VEM}}(\cdot,\cdot)$ are only element-local, i.e. we can again apply static condensation, remove interior degrees of freedom and obtain a global system for the vertex and edge unknowns in 2D or vertex, edge and face unknowns in 3D.
In contrast to the DG methods and the HDG methods unknowns are now associated with different geometrical entities (not just elements or facets).
In \cref{fig:VEMcouplings} the couplings of the VEM unknowns (after static condensation) are illustrated.

\section{Sparsity measures and basic coupling relations} \label{sec:measures}
The classical framework of the Laplace problem with Dirichlet boundary conditions in \cref{sec:model_problem} is excellent for presenting the different methods in their simplest form.
However, it adds additional complexity of the discussion of boundary effects when comparing the coupling relations of the different methods in \cref{sec:measures}. 
To avoid overcomplication we will therefore consider \emph{periodic boundaries} instead so that all geometrical entities and d.o.f.s. are \emph{interior} entities. Moreover, we consider simple periodic meshes where the unit cell (square or cube) is filled by a small number of elements.

\subsection{Notation}
To compare the computational costs of the different methods, we introduce the following notation.

We require some notation for geometrical entities of the mesh and a corresponding neighborhood relation. 
For elements ($d$ dimensional) and facets ($d-1$ dimensional), we use the notation El and Ft whereas for vertices (0 dimensional), edges (1 dimensional), faces (2 dimensional) and cells (3 dimensional) we use the notation V, Ed, Fa and C.

$N_{S}$ denotes the number of entities in the mesh with the same dimension, with $S \in \mathcal{S} = \{\text{V},\text{Ed},\text{Fa}, \text{C}\}$ or $S \in \{\text{El}, \text{Ft}\}$. Especially, by $\Nel$ we denote the number of elements in a mesh and by $\Nft$ the number of facets in the mesh.
An obvious measure for the computational costs is the number of degrees of freedom (\ndof) which we denote by $\ndof$. More important in view of the sparsity of the resulting linear system are the following two quantities:
\begin{align*}
    \ncdof & \quad \text{represents the number of \emph{coupling} degrees of freedom that remain after static condensation, and } \\
    \nnze & \quad \text{signifies the number of non-zero entries in the linear system after static condensation.  } %
\end{align*}
We point out that the $\nnze$ we consider here is the maximum that is reached if all coupling degrees of freedom produce a non-zero entry. 
In practice, the number of non-zero entries can be smaller, e.g. if one uses specific orthogonal basis functions.

\paragraph{Geometrical entities with same local topology.}
In preparation for the discussion of unknowns and couplings in the HDG and VE methods -- which have the most complex coupling structure -- we introduce a generic approach for categorizing geometric entities with the same local topology.
We denote by $\mathcal{E}$ the set of the types of geometrical entities \emph{with the same local topology}. Each element in $\mathcal{E}$ is a pair of a geometrical entity of fixed dimension $E \in \mathcal{S}$ and an integer, e.g. $(\text{Ed},1)$ would denote the first type of edges in the mesh, $(\text{Fa},2)$ the second type of faces, etc. 
In the simplest case, we can have $\mathcal{E} = \{ (\text{V},1), (\text{Ed},1), (\text{Fa},1) \}$ for a periodic quadrilateral mesh in 2D, cf. \cref{sec:quadrilaterals}, where all vertices, edges and faces have the same topology. 
Entities of the same dimension can also appear several times with different local topologies. For instance the truncated octahedron below, cf. \cref{sec:truncatedoctahedron}, has faces with 4 edges and faces with 6 edges. In that case we have $\mathcal{E} = \{ (\text{V},1), (\text{Ed},1), (\text{Fa},1), (\text{Fa},2), (\text{C},1) \}$. 
For entities of a specified dimension we introduce the notation $\mathcal{E}_m$ for the set of integers of $\mathcal{E}$ corresponding to geometrical entities of dimension $m$ and set $\mathcal{E}_{\text{El}} = \mathcal{E}_d$ and $\mathcal{E}_{\text{Ft}} = \mathcal{E}_{d-1}$.

\paragraph{Ratio between entities of a specific type and elements.}
We will normalize the measures -- especially \nnze~ -- by the number of elements in the mesh. On different meshes, the number of certain geometrical entities per element can vary. 
With $R_{\text{El}}^P$ we denote the ratio between the number of elements of type $P \in \mathcal{E}$ and the number of elements in the mesh. 

\paragraph{Neighborhood relation between mesh entities.}
We say that a geometrical entity is a neighbor of another entity if they are part of the same element in the mesh. 
We denote the number of neighboring entities of fixed dimension $P \in \mathcal{S}$ to a geometrical entity of type topology type $Q$ in the mesh by $\Nb{Q}{P}$. For instance, $\Nb{(\text{Fa},2)}{\text{Ed}}$ describes how many edges are neighboring each face in the second type of faces and $\Nb{(\text{El},1)}{\text{Ft}}$ how many facets are attached to each element in the first type of elements.

For the different methods, we will focus on comparing the \ncdof, which are unknowns entering the global linear system (possibly after static condensation) and are therefore relevant for the analysis of computational costs. 
Further, we consider the non-zero entries $\nnze$ related to the \ncdof, i.e. after condensation if applicable. 
They are a crucial metric for assessing computational efficiency. 

In general, for \ndof, \ncdof, and \nnze a subscript will indicate the method, e.g. $\ndof_{\text{DG}}$ for the DG method.
To express the global \ncdof and \nnze, we will rely on local (per geometrical entity) degrees of freedom, which we denote by $\ndof^{\text{El}}$, $\ndof^{\text{Ft}}$, $\ndof^{\text{V}}$, $\ndof^{\text{Ed}}$ or $\ndof^{\text{Fa}}$. 

\subsection{Number of coupling degrees of freedom (\ncdof)}
\begin{subequations}\label{eq:ncdofs}
  For DG and the Trefftz DG methods, the global \ncdof are simply local $\ndof^{\text{El}}$ summed up over all elements. 
\begin{align}
  \ncdof_{\text{DGM}} &= \ndof_{\text{DGM}} = \Nel \cdot \ndof_{\text{DGM}}^{\text{El}}, &\quad \text{DGM} \in &\{ \text{DG}, \text{TDG1}, \text{TDG2} \}
  \intertext{
  For Hybrid DG and HHO methods the volume degrees of freedom can be condensed out, so that the global \ncdof are the local $\ndof^{\text{Ft}}$ summped up over the facets. 
  }
	\ncdof_{\text{HDGM}} &= \Nft \cdot \ndof_{\text{HDGM}}^{\text{Ft}}, &\quad \text{HDGM} \in &\{ \text{HDG}, \text{HHO} \}
  \intertext{
  For VEM, we also condense the volume degrees of freedom, leaving the degrees of freedom on vertices and edges (and faces in 3D).
  With this generic notation for different types of geometrical entities introduced above, we can write the global \ncdof for the VEM method as \vspace*{-0.5cm}
  }
  \ncdof_{\text{VEM}} &= \overbrace{\sum_{(X,i)\in\mathcal{E}} N_{(X,i)} \cdot \ndof_{\text{VEM}}^{X}}^{ = \ndof_{\text{VEM}}} - \Nel \cdot \ndof_{\text{VEM}}^{\text{El}}.
\end{align}
\end{subequations}
In \cref{tab:ndof} we summarize the local \ncdof for the different methods and (relevant) geometrical entities. We note that in the considered methods the number of unknowns per geometrical entity only depends on the dimension of the entity and not on the type of the entity (if different types exist).

  \begin{table}[!ht]
	\begin{center}
	\begin{tabular}{lllll}
	  \toprule
	  method &  abbr. & $d$ dim. & 2D & 3D\\ 
	 \midrule
	 DG & $\ncdof_{\text{DG}}^{\text{El}}$ & $ \begin{psmallmatrix} k + d \\ d \end{psmallmatrix} $ & $\frac{(k+1)(k+2)}{2}$ & $\frac{(k+1)(k+2)(k+3)}{6}$ \\[1ex]
	 \rowcolor{gray!20} Trefftz DG 2 &  $\ncdof_{\text{TDG2}}^{\text{El}}$ & $ \begin{psmallmatrix} k + d  \\ d \end{psmallmatrix}  -  \begin{psmallmatrix} k-2+d \\ d \end{psmallmatrix}  $ & $2k+1$ & $(k+1)^2$ \\ ~ \\[-1.8ex]
	 Trefftz DG 1 & $\ncdof_{\text{TDG1}}^{\text{El}}$ &  $ \begin{psmallmatrix} k + d  \\ d \end{psmallmatrix}  -  \begin{psmallmatrix} k-1+d \\ d \end{psmallmatrix}  $ & $k+1$ & $\frac{(k+1)(k+2)}{2}$ \\
	 \midrule 
	 \rowcolor{gray!20} Hybrid DG & $\ncdof_{\text{HDG}}^{\text{Ft}}$  &  $ \begin{psmallmatrix} k + d - 1 \\ d-1 \end{psmallmatrix} $ & $k+1$ & $\frac{(k+1)(k+2)}{2}$ \\~ \\[-1.8ex]
	 Hybrid HO &  $\ncdof_{\text{HHO}}^{\text{Ft}}$  &$ \begin{psmallmatrix} k -1 + d - 1 \\ d-1 \end{psmallmatrix} $ & $k$ & $\frac{k(k+1)}{2}$ \\ ~ \\[-1.8ex]
	 \rowcolor{gray!20} VEM & $\ncdof_{\text{VEM}}^{\text{V}}$ & $1$ & $1$ & $1$ \\[1ex]
	 \rowcolor{gray!20} & $\ncdof_{\text{VEM}}^{\text{Ed}}$ & $k-1$ & $k-1$ & $k-1$ \\[1ex]
	 \rowcolor{gray!20} & $\ncdof_{\text{VEM}}^{\text{Fa}}$ & $\frac{(k-1)k}{2}$ & -- & $\frac{(k-1)k}{2}$  \\[1ex]
	 \rowcolor{gray!20} & $\ncdof_{\text{VEM}}^{\text{C}}$ & -- & -- & --  \\[1ex]
	 \bottomrule
	\end{tabular}
	\end{center}\vspace*{-0.5cm}
	\caption{Number of coupling degrees of freedom for the considered methods on one element, facet, vertex, edge or face, respectively.}
	\label{tab:ndof}
	\end{table}

  \subsection{Number of non-zero entries in the system matrix (\nnze)}
  The \nnze for the DG and Trefftz DG method follows the simple pattern that all degrees of freedom of one element couple with all degrees of freedom of neighboring elements. This yields
  \begin{subequations} \label{eq:nnzes}
  \begin{align}
  \nnze_{\text{DGM}}/N_{\text{El}} &= \sum_{i \in \mathcal{E}_{\text{El}}} R_{\text{El}}^{(\text{El},i)} (\Nb{(\text{El},i)}{\text{Ft}}+1) (\ndof_{\text{DGM}}^{\text{El}})^2, \quad \text{DGM} \in \{ \text{DG}, \text{TDG1}, \text{TDG2} \}.\\
  \intertext{For the HDG methods the \nnze are determined by the couplings of the facet unknowns to neighboring facet unknowns. We can however have types of facets with different numbers of neighbors within the mesh.}
  \nnze_{\text{HDGM}}/N_{\text{El}} &= \sum_{i \in \mathcal{E}_{\text{Ft}}} R_{\text{El}}^{(\text{Ft},i)} \Nb{(\text{Ft},i)}{\text{Ft}} (\ndof_{\text{HDGM}}^{\text{Ft}})^2, \quad \text{HDGM} \in \{ \text{HDG}, \text{HHO} \}.
\intertext{ For the Virtual Element Method, the formula for the \nnze is a bit more complex. Each group of unknowns associated with a geometrical entity couples with all groups of unknowns associated with neighboring (non-volumetric) geometrical entities. }
  \nnze_{\text{VEM}}/N_{\text{El}} &= \sum_{(X,i) \in \mathcal{E}} R_{\text{El}}^{(X,i)} \cdot \ncdof_{\text{VEM}}^{X} \sum_{S \in \mathcal{S}} \Nb{(X,i)}{S} \cdot \ncdof_{\text{VEM}}^{S}
  \end{align}
\end{subequations}

\section{Sparsity comparison on periodic polytopal meshes} \label{sec:periodicpolys}
In this section, we compare the sparsity of the linear system arising from the different methods.
Note that all methods above show the same optimal convergence rates (in suitable norms) for several problems, including the model problem of \cref{sec:model_problem}, with an error that only differs by a constant. Therefore comparison of the quantities of $\ncdof,\nnze$ for fixed $k$ is sensible.
We compare all six previously introduced methods, but mark the HHO in gray as it is only applicable in special situations, as well as the Trefftz DG1 method, as second order problems are more common.

\newpage

\subsection{2D: Triangles}
\begin{figure}[ht!]
	\begin{center} 
		\begin{minipage}{0.48\textwidth}
			\centering
		\begin{tikzpicture}[scale=1.5]
			\foreach \x in {0,0.33,0.66,0.99,1.32,1.65,1.98,2.31,2.64}{
				\foreach \y in {0,0.33,0.66,0.99,1.32,1.65,1.98,2.31,2.64}{
					\pgfmathparse{35*rnd+65}
					\edef\col{\pgfmathresult}
					\fill[white!\col!black] (\x,\y) -- (\x+0.33,\y) -- (\x,\y+0.33) -- cycle;
					\pgfmathparse{35*rnd+65}
					\edef\col{\pgfmathresult}
					\fill[white!\col!black] (\x+0.33,\y) -- (\x+0.33,\y+0.33) -- (\x,\y+0.33) -- cycle;
				}
			}
			\def\nonperiodic{0.15mm}
		  \foreach \start in {0,0.33,0.66,0.99,1.32,1.65,1.98,2.31,2.64,2.97}{
			  \draw[line width=\nonperiodic] (0,\start) -- (2.97,\start); %
			  \draw[line width=\nonperiodic] (\start, 0) -- (\start, 2.97); %
		  }
		  \foreach \val in {0.33,0.66,0.99,1.32,1.65,1.98,2.31,2.64,2.97}{
			  \draw[line width=\nonperiodic] (0,\val) -- (\val, 0);
		  }
		  \foreach \val in {0.33,0.66,0.99,1.32,1.65,1.98,2.31,2.64}{
			  \draw[line width=\nonperiodic] (\val, 2.97) -- (2.97, \val);
			}
		\end{tikzpicture}
		\end{minipage}
			\begin{minipage}{0.48\textwidth}
				\centering
				\begin{tabular}{rcccc}
					& & $\Nb{(X,i)}{\text{Ft}}$ &  & \\
					& & \rotatebox{90}{$=$} &  & \\
					\toprule
					$(X,i) \downarrow $  & $\Nb{(X,i)}{\text{V}}$ & $\Nb{(X,i)}{\text{Ed}}$ & \color{gray} $\Nb{(X,i)}{\text{Fa}}$ & $R^{(X,i)}_{\text{El}}$\\
					\midrule
					(V,1) & 7 & 12 & \color{gray}6 & $\frac12$
					\\[0.5ex]
					(Ft,1) = (Ed,1) & 4 & 5 & \color{gray}2 & $\frac32$
					\\[0.5ex]
					(El,1) = (Fa,1) & \color{gray}3 &  3 &  \color{gray}1 & \color{gray}1
					\\
					\bottomrule 
				\end{tabular}
			\end{minipage}
	\end{center}	\vspace*{-0.5cm}
	\caption{Sketch of structured periodic triangle mesh (left) and the relevant neighborhood topology numbers (right).}
	\label{fig:triangle}
  \end{figure}  
We consider a periodic triangle mesh that is obtained by dividing the periodic square across one diagonal, cf. \cref{fig:triangle} for a sketch. We only have one type of vertices, edges and faces, respectively, i.e. $\{1\} = \mathcal{E}_0 = \mathcal{E}_1 = \mathcal{E}_2$. We can compute the number of coupling unknowns and non-zero entries in the linear system per element for the different methods with \eqref{eq:ncdofs} and \eqref{eq:nnzes} and the topology numbers provided in \cref{fig:triangle} which yields the result in \cref{tab:triangle:nnnpol}.
\begin{table}[ht!]
\begin{center}	
\begin{NiceTabular}{rr@{~}r@{}r@{~}r@{}r|r@{~}r@{}r@{~}r@{}r@{~}r@{}r@{~}r@{}r}[colortbl-like]
\toprule
& \multicolumn{5}{c}{$\texttt{ncdof}/N_{El}$}& \multicolumn{9}{c}{$\texttt{nnze}/N_{El}$} \\
\midrule
DG$ $&$  \frac{1}{2} k^2$&$  + $&$  \frac{3}{2} k  $&$  + $&$ 1$&$  k^4$&$  + $&$ 6 k^3$&$  + $&$ 13 k^2$&$  + $&$ 12 k  $&$  + $&$ 4$\\[0.5ex]\rowcolor{gray!20}\\[-2ex]
\rowcolor{gray!20}TDG2$ $&$ $&$ $&$ 2 k  $&$  + $&$ 1$&$ $&$ $&$ $&$ $&$ 16 k^2$&$  + $&$ 16 k  $&$  + $&$ 4$\\[0.5ex]\\[-2ex]
TDG1$ $&$ $&$ $&$  k  $&$  + $&$ 1$&$ $&$ $&$ $&$ $&$ 4 k^2$&$  + $&$ 8 k  $&$  + $&$ 4$\\[0.5ex]\rowcolor{gray!20}\\[-2ex]
\rowcolor{gray!20}HDG$ $&$ $&$ $&$  \frac{3}{2} k  $&$  + $&$  \frac{3}{2}$&$ $&$ $&$ $&$ $&$  \frac{15}{2} k^2$&$  + $&$ 15 k  $&$  + $&$  \frac{15}{2}$\\[0.5ex]\\[-2ex]
HHO$ $&$ $&$ $&$  \frac{3}{2} k  $&$ $&$ $&$ $&$ $&$ $&$ $&$  \frac{15}{2} k^2$&$ $&$ $&$ $&$ $\\[0.5ex]\rowcolor{gray!20}\\[-2ex]
\rowcolor{gray!20}VEM$ $&$ $&$ $&$  \frac{3}{2} k  $&$  - $&$ 1$&$ $&$ $&$ $&$ $&$  \frac{15}{2} k^2$&$  - $&$ 3 k  $&$  - $&$ 1$\\\bottomrule
\end{NiceTabular}

\end{center}\vspace*{-0.5cm}
\caption{Number of coupling unknowns and non-zero entries per element for different methods on a periodic triangle mesh.} 
\label{tab:triangle:nnnpol}
\end{table}

We observe that all methods, except the standard DG approach, have the number of non-zero entries scaling like $\mathcal{O}(k^2)$, while the DG method has a higher scaling of $\mathcal{O}(k^4)$. The constant factor in front of the $k^2$ term is different for the different methods and strongly in favor of the methods that condense out the volume degrees of freedom, namely the HDG, HHO and VEM methods -- unless we have a first order operator and can apply the TDG1 method.
In \cref{tab:triangle:nnn} we list the concrete numbers of coupling dofs and non-zero entries in the linear system per element for the different methods on the periodic triangle mesh for $k=1,\ldots,10$.
\begin{table}[!ht]
\centering
\begin{NiceTabular}{R{2.5cm}@{}R{1.05cm}@{}R{1.1cm}@{}R{1.15cm}@{}R{1.2cm}@{}R{1.25cm}@{}R{1.3cm}@{}R{1.35cm}@{}R{1.4cm}@{}R{1.45cm}@{}R{1.5cm}@{}}[colortbl-like]
    \toprule
\text{method}$~\downarrow$ $\setminus~k \rightarrow$ & \multicolumn{1}{r}{1}& \multicolumn{1}{r}{2}& \multicolumn{1}{r}{3}& \multicolumn{1}{r}{4}& \multicolumn{1}{r}{5}& \multicolumn{1}{r}{6}& \multicolumn{1}{r}{7}& \multicolumn{1}{r}{8}& \multicolumn{1}{r}{9}& \multicolumn{1}{r}{10}\\
   \midrule
$\texttt{ndof}_{\text{DG}}/N_{\text{el}}$&         3\hphantom{.0} &          6\hphantom{.0} &         10\hphantom{.0} &         15\hphantom{.0} &         21\hphantom{.0} &         28\hphantom{.0} &         36\hphantom{.0} &         45\hphantom{.0} &         55\hphantom{.0} &         66\hphantom{.0}\\
\rowcolor{gray!20}$\texttt{ndof}_{\text{TDG2}}/N_{\text{el}}$&         3\hphantom{.0} &          5\hphantom{.0} &          7\hphantom{.0} &          9\hphantom{.0} &         11\hphantom{.0} &         13\hphantom{.0} &         15\hphantom{.0} &         17\hphantom{.0} &         19\hphantom{.0} &         21\hphantom{.0}\\
\color{gray}$\texttt{ndof}_{\text{TDG1}}/N_{\text{el}}$&\color{gray}         2\hphantom{.0} & \color{gray}         3\hphantom{.0} & \color{gray}         4\hphantom{.0} & \color{gray}         5\hphantom{.0} & \color{gray}         6\hphantom{.0} & \color{gray}         7\hphantom{.0} & \color{gray}         8\hphantom{.0} & \color{gray}         9\hphantom{.0} & \color{gray}        10\hphantom{.0} & \color{gray}        11\hphantom{.0}\\
\rowcolor{gray!20}$\texttt{ncdof}_{\text{HDG}}/N_{\text{el}}$&         3\hphantom{.0} &        4.5 &          6\hphantom{.0} &        7.5 &          9\hphantom{.0} &       10.5 &         12\hphantom{.0} &       13.5 &         15\hphantom{.0} &       16.5\\
\color{gray}$\texttt{ncdof}_{\text{HHO}}/N_{\text{el}}$&\color{gray}       1.5 & \color{gray}         3\hphantom{.0} & \color{gray}       4.5 & \color{gray}         6\hphantom{.0} & \color{gray}       7.5 & \color{gray}         9\hphantom{.0} & \color{gray}      10.5 & \color{gray}        12\hphantom{.0} & \color{gray}      13.5 & \color{gray}        15\hphantom{.0}\\
\rowcolor{gray!20}$\texttt{ncdof}_{\text{VEM}}/N_{\text{el}}$&       0.5 &          2\hphantom{.0} &        3.5 &          5\hphantom{.0} &        6.5 &          8\hphantom{.0} &        9.5 &         11\hphantom{.0} &       12.5 &         14\hphantom{.0}\\
    \bottomrule
\end{NiceTabular}

\begin{NiceTabular}{R{2.5cm}@{}R{1.05cm}@{}R{1.1cm}@{}R{1.15cm}@{}R{1.2cm}@{}R{1.25cm}@{}R{1.3cm}@{}R{1.35cm}@{}R{1.4cm}@{}R{1.45cm}@{}R{1.5cm}@{}}[colortbl-like]
    \toprule
\text{method}$~\downarrow$ $\setminus~k \rightarrow$ & \multicolumn{1}{r}{1}& \multicolumn{1}{r}{2}& \multicolumn{1}{r}{3}& \multicolumn{1}{r}{4}& \multicolumn{1}{r}{5}& \multicolumn{1}{r}{6}& \multicolumn{1}{r}{7}& \multicolumn{1}{r}{8}& \multicolumn{1}{r}{9}& \multicolumn{1}{r}{10}\\
   \midrule
$\texttt{nnze}_{\text{DG}}/N_{\text{el}}$&        36\hphantom{.0} &        144\hphantom{.0} &        400\hphantom{.0} &        900\hphantom{.0} &       1764\hphantom{.0} &       3136\hphantom{.0} &       5184\hphantom{.0} &       8100\hphantom{.0} &      12100\hphantom{.0} &      17424\hphantom{.0}\\
\rowcolor{gray!20}$\texttt{nnze}_{\text{TDG2}}/N_{\text{el}}$&        36\hphantom{.0} &        100\hphantom{.0} &        196\hphantom{.0} &        324\hphantom{.0} &        484\hphantom{.0} &        676\hphantom{.0} &        900\hphantom{.0} &       1156\hphantom{.0} &       1444\hphantom{.0} &       1764\hphantom{.0}\\
\color{gray}$\texttt{nnze}_{\text{TDG1}}/N_{\text{el}}$&\color{gray}        16\hphantom{.0} & \color{gray}        36\hphantom{.0} & \color{gray}        64\hphantom{.0} & \color{gray}       100\hphantom{.0} & \color{gray}       144\hphantom{.0} & \color{gray}       196\hphantom{.0} & \color{gray}       256\hphantom{.0} & \color{gray}       324\hphantom{.0} & \color{gray}       400\hphantom{.0} & \color{gray}       484\hphantom{.0}\\
\rowcolor{gray!20}$\texttt{nnze}_{\text{HDG}}/N_{\text{el}}$&        30\hphantom{.0} &       67.5 &        120\hphantom{.0} &      187.5 &        270\hphantom{.0} &      367.5 &        480\hphantom{.0} &      607.5 &        750\hphantom{.0} &      907.5\\
\color{gray}$\texttt{nnze}_{\text{HHO}}/N_{\text{el}}$&\color{gray}       7.5 & \color{gray}        30\hphantom{.0} & \color{gray}      67.5 & \color{gray}       120\hphantom{.0} & \color{gray}     187.5 & \color{gray}       270\hphantom{.0} & \color{gray}     367.5 & \color{gray}       480\hphantom{.0} & \color{gray}     607.5 & \color{gray}       750\hphantom{.0}\\
\rowcolor{gray!20}$\texttt{nnze}_{\text{VEM}}/N_{\text{el}}$&       3.5 &         23\hphantom{.0} &       57.5 &        107\hphantom{.0} &      171.5 &        251\hphantom{.0} &      345.5 &        455\hphantom{.0} &      579.5 &        719\hphantom{.0}\\
    \bottomrule
\end{NiceTabular}

\caption{\ncdof and \nnze per element for different methods on periodic triangle mesh for
$k\!\!=\!\!1,\!..,\!10$ (rounded up to one decimal place).}
\label{tab:triangle:nnn}
\end{table}

\newpage
\subsection{2D: Quadrilaterals} \label{sec:quadrilaterals}
\begin{figure}[ht!]
	\begin{center} 
		\begin{minipage}{0.48\textwidth}
			\centering
	\begin{tikzpicture}[scale=1.5]
		  \def\nonperiodic{0.15mm}
		  \foreach \x in {0,0.33,0.66,0.99,1.32,1.65,1.98,2.31,2.64}{
			  \foreach \y in {0,0.33,0.66,0.99,1.32,1.65,1.98,2.31,2.64}{
				  \pgfmathparse{35*rnd+65}
				  \edef\col{\pgfmathresult}
				  \fill[white!\col!black] (\x,\y) -- (\x+0.33,\y) -- (\x+0.33,\y+0.33) -- (\x,\y+0.33) -- cycle;
			  }
		  }
		  \foreach \start in {0,0.33,0.66,0.99,1.32,1.65,1.98,2.31,2.64,2.97}{
			  \draw[line width=\nonperiodic] (0,\start) -- (2.97,\start); %
			  \draw[line width=\nonperiodic] (\start, 0) -- (\start, 2.97); %
		  }
	\end{tikzpicture}  

\end{minipage}
\begin{minipage}{0.48\textwidth}
	\centering
\begin{tabular}{rcccc}
		& & $\Nb{(X,i)}{\text{Ft}}$ & & \\
		& & \rotatebox{90}{$=$} & & \\
		\toprule
		$(X,i) \downarrow $  & $\Nb{(X,i)}{\text{V}}$ & $\Nb{(X,i)}{\text{Ed}}$ & \color{gray} $\Nb{(X,i)}{\text{Fa}}$ & $R^{(X,i)}_{\text{El}}$\\
		\midrule
		(V,1) & 9 & 12 & \color{gray}4 & 1 
		\\ 
		(Ft,1) = (Ed,1) & 6 & 7 & \color{gray}2 & 2 
		\\
		(El,1) = (Fa,1) & \color{gray}4 & 4 &  \color{gray}1 & \color{gray}1
		\\
		\bottomrule 
	\end{tabular}
\end{minipage}
\end{center}\vspace*{-0.5cm}
\caption{Sketch of structured periodic quadrilateral mesh (left) and the relevant neighborhood topology numbers (right).}
	\label{fig:quads}
\end{figure}  
For the periodic quadrilateral mesh, cf. \cref{fig:quads}, we also only have one type of vertices, edges and faces, respectively, i.e. $\{1\} = \mathcal{E}_0 = \mathcal{E}_1 = \mathcal{E}_2$. We obtain the dependency on $k$ in \ncdof and \nnze as shown in \cref{tab:quadrilateral:nnnpol}.
\begin{table}[ht!]
\begin{center}	
\begin{NiceTabular}{rr@{~}r@{}r@{~}r@{}r|r@{~}r@{}r@{~}r@{}r@{~}r@{}r@{~}r@{}r}[colortbl-like]
\toprule
& \multicolumn{5}{c}{$\texttt{ncdof}/N_{El}$}& \multicolumn{9}{c}{$\texttt{nnze}/N_{El}$} \\
\midrule
DG$ $&$  \frac{1}{2} k^2$&$  + $&$  \frac{3}{2} k  $&$  + $&$ 1$&$  \frac{5}{4} k^4$&$  + $&$  \frac{15}{2} k^3$&$  + $&$  \frac{65}{4} k^2$&$  + $&$ 15 k  $&$  + $&$ 5$\\[0.5ex]\rowcolor{gray!20}\\[-2ex]
\rowcolor{gray!20}TDG2$ $&$ $&$ $&$ 2 k  $&$  + $&$ 1$&$ $&$ $&$ $&$ $&$ 20 k^2$&$  + $&$ 20 k  $&$  + $&$ 5$\\[0.5ex]\\[-2ex]
TDG1$ $&$ $&$ $&$  k  $&$  + $&$ 1$&$ $&$ $&$ $&$ $&$ 5 k^2$&$  + $&$ 10 k  $&$  + $&$ 5$\\[0.5ex]\rowcolor{gray!20}\\[-2ex]
\rowcolor{gray!20}HDG$ $&$ $&$ $&$ 2 k  $&$  + $&$ 2$&$ $&$ $&$ $&$ $&$ 14 k^2$&$  + $&$ 28 k  $&$  + $&$ 14$\\[0.5ex]\\[-2ex]
HHO$ $&$ $&$ $&$ 2 k  $&$ $&$ $&$ $&$ $&$ $&$ $&$ 14 k^2$&$ $&$ $&$ $&$ $\\[0.5ex]\rowcolor{gray!20}\\[-2ex]
\rowcolor{gray!20}VEM$ $&$ $&$ $&$ 2 k  $&$  - $&$ 1$&$ $&$ $&$ $&$ $&$ 14 k^2$&$  - $&$ 4 k  $&$  - $&$ 1$\\\bottomrule
\end{NiceTabular}

\end{center}\vspace*{-0.5cm}
\caption{Number of coupling unknowns and non-zero entries per element for different methods on a periodic quadrilateral mesh.} 
\label{tab:quadrilateral:nnnpol}
\end{table}

We observe the same scaling as for the triangle case where the constants in front of the $k^2$ term in the \nnze are closer to each other for the different methods. Moreover, in \cref{tab:quadrilateral:nnn} we observe that especially for lower orders the difference becomes less significant and even slightly in favor of the Trefftz DG method.
  
\begin{table}[!ht]
    \centering
	\begin{NiceTabular}{R{2.5cm}@{}R{1.05cm}@{}R{1.1cm}@{}R{1.15cm}@{}R{1.2cm}@{}R{1.25cm}@{}R{1.3cm}@{}R{1.35cm}@{}R{1.4cm}@{}R{1.45cm}@{}R{1.5cm}@{}}[colortbl-like]
    \toprule
\text{method}$~\downarrow$ $\setminus~k \rightarrow$ & \multicolumn{1}{r}{1}& \multicolumn{1}{r}{2}& \multicolumn{1}{r}{3}& \multicolumn{1}{r}{4}& \multicolumn{1}{r}{5}& \multicolumn{1}{r}{6}& \multicolumn{1}{r}{7}& \multicolumn{1}{r}{8}& \multicolumn{1}{r}{9}& \multicolumn{1}{r}{10}\\
   \midrule
$\texttt{ndof}_{\text{DG}}/N_{\text{el}}$&         3 &          6 &         10 &         15 &         21 &         28 &         36 &         45 &         55 &         66\\
\rowcolor{gray!20}$\texttt{ndof}_{\text{TDG2}}/N_{\text{el}}$&         3 &          5 &          7 &          9 &         11 &         13 &         15 &         17 &         19 &         21\\
\color{gray}$\texttt{ndof}_{\text{TDG1}}/N_{\text{el}}$&\color{gray}         2 & \color{gray}         3 & \color{gray}         4 & \color{gray}         5 & \color{gray}         6 & \color{gray}         7 & \color{gray}         8 & \color{gray}         9 & \color{gray}        10 & \color{gray}        11\\
\rowcolor{gray!20}$\texttt{ncdof}_{\text{HDG}}/N_{\text{el}}$&         4 &          6 &          8 &         10 &         12 &         14 &         16 &         18 &         20 &         22\\
\color{gray}$\texttt{ncdof}_{\text{HHO}}/N_{\text{el}}$&\color{gray}         2 & \color{gray}         4 & \color{gray}         6 & \color{gray}         8 & \color{gray}        10 & \color{gray}        12 & \color{gray}        14 & \color{gray}        16 & \color{gray}        18 & \color{gray}        20\\
\rowcolor{gray!20}$\texttt{ncdof}_{\text{VEM}}/N_{\text{el}}$&         1 &          3 &          5 &          7 &          9 &         11 &         13 &         15 &         17 &         19\\
    \bottomrule
\end{NiceTabular}

	\begin{NiceTabular}{R{2.5cm}@{}R{1.05cm}@{}R{1.1cm}@{}R{1.15cm}@{}R{1.2cm}@{}R{1.25cm}@{}R{1.3cm}@{}R{1.35cm}@{}R{1.4cm}@{}R{1.45cm}@{}R{1.5cm}@{}}[colortbl-like]
    \toprule
\text{method}$~\downarrow$ $\setminus~k \rightarrow$ & \multicolumn{1}{r}{1}& \multicolumn{1}{r}{2}& \multicolumn{1}{r}{3}& \multicolumn{1}{r}{4}& \multicolumn{1}{r}{5}& \multicolumn{1}{r}{6}& \multicolumn{1}{r}{7}& \multicolumn{1}{r}{8}& \multicolumn{1}{r}{9}& \multicolumn{1}{r}{10}\\
   \midrule
$\texttt{nnze}_{\text{DG}}/N_{\text{el}}$&        45 &        180 &        500 &       1125 &       2205 &       3920 &       6480 &      10125 &      15125 &      21780\\
\rowcolor{gray!20}$\texttt{nnze}_{\text{TDG2}}/N_{\text{el}}$&        45 &        125 &        245 &        405 &        605 &        845 &       1125 &       1445 &       1805 &       2205\\
\color{gray}$\texttt{nnze}_{\text{TDG1}}/N_{\text{el}}$&\color{gray}        20 & \color{gray}        45 & \color{gray}        80 & \color{gray}       125 & \color{gray}       180 & \color{gray}       245 & \color{gray}       320 & \color{gray}       405 & \color{gray}       500 & \color{gray}       605\\
\rowcolor{gray!20}$\texttt{nnze}_{\text{HDG}}/N_{\text{el}}$&        56 &        126 &        224 &        350 &        504 &        686 &        896 &       1134 &       1400 &       1694\\
\color{gray}$\texttt{nnze}_{\text{HHO}}/N_{\text{el}}$&\color{gray}        14 & \color{gray}        56 & \color{gray}       126 & \color{gray}       224 & \color{gray}       350 & \color{gray}       504 & \color{gray}       686 & \color{gray}       896 & \color{gray}      1134 & \color{gray}      1400\\
\rowcolor{gray!20}$\texttt{nnze}_{\text{VEM}}/N_{\text{el}}$&         9 &         47 &        113 &        207 &        329 &        479 &        657 &        863 &       1097 &       1359\\
    \bottomrule
\end{NiceTabular}

	\caption{\ncdof and \nnze per element for different methods on periodic quadrilateral mesh for $k=1,\ldots,10$.}
	\label{tab:quadrilateral:nnn}
\end{table}

\newpage	
\subsection{2D: Hexagons}

\begin{figure}[ht!]
	\begin{center} 
		\begin{minipage}{0.48\textwidth}
			\centering
	\begin{tikzpicture}[scale=1.5]
			\def\periodic{0.12mm}
			\draw[line width=\periodic, color=gray] (0,0) -- (2.97,0);
			\draw[line width=\periodic, color=gray] (0,2.97) -- (2.97,2.97);
			\draw[line width=\periodic, color=gray] (0,0) -- (0,2.97);
			\draw[line width=\periodic, color=gray] (2.97,0) -- (2.97,2.97);

		  \def\nonperiodic{0.15mm}
		  \fill[red!40] (0,0) -- (0.33,0) -- (0,0.33) -- cycle;
		  \fill[red!40] (2.97,2.97) -- (2.97,2.64) -- (2.64,2.97) -- cycle;
		  \fill[red!40] (2.97,0) -- (2.97,0.33) -- (2.64,0.33) -- (2.64,0) -- cycle;
		  \fill[red!40] (0,2.97) -- (0.33,2.97) -- (0.33,2.64) -- (0,2.64) -- cycle;
		  
		  \fill[green!40] (0.66,0) -- (1.32,0) -- (0.99,0.33) -- (0.66,0.33) -- cycle;
		  \fill[green!40] (0.66,2.97) -- (1.32,2.97) -- (1.32,2.64) -- (0.99,2.64) -- cycle;
		  
		  \fill[blue!40] (1.65,0) -- (2.31,0) -- (1.98,0.33) -- (1.65,0.33) -- cycle;
		  \fill[blue!40] (1.65,2.97) -- (2.31,2.97) -- (2.31,2.64) -- (1.98,2.64) -- cycle;
		  
		  \fill[cyan!40] (0,0.66) -- (0.33,0.66) -- (0.33,0.99) -- (0,1.32) -- cycle;
		  \fill[cyan!40] (2.97,0.66) -- (2.64,0.99) -- (2.64,1.32) -- (2.97,1.32) -- cycle;
		  
		  \fill[magenta!40] (0,1.65) -- (0.33,1.65) -- (0.33,1.98) -- (0,2.31) -- cycle;
		  \fill[magenta!40] (2.97,1.65) -- (2.64,1.98) -- (2.64,2.31) -- (2.97,2.31) -- cycle;

		  \def\nonperiodic{0.225mm}
		  \foreach \valtmp in {0,0.99,1.98}{
			  \foreach \tmp in {0.33,1.32,2.31}{
				  \draw[line width=\nonperiodic] (\tmp,\valtmp) -- (\tmp+0.33,\valtmp) -- (\tmp+0.33,\valtmp+0.33) -- (\tmp,\valtmp+0.66) -- (\tmp-0.33,\valtmp+0.66) -- (\tmp-0.33,\valtmp+0.33) -- cycle;
				  \pgfmathparse{50*rnd+50}
				  \edef\col{\pgfmathresult}
				  \fill[white!\col!black] (\tmp,\valtmp) -- (\tmp+0.33,\valtmp) -- (\tmp+0.33,\valtmp+0.33) -- (\tmp,\valtmp+0.66) -- (\tmp-0.33,\valtmp+0.66) -- (\tmp-0.33,\valtmp+0.33) -- cycle;
			  }
			  \foreach \tmp in {0.66,1.65,2.64}{
				  \draw[line width=\nonperiodic] (\tmp,\valtmp+0.33) -- (\tmp+0.33,\valtmp+0.33) -- (\tmp+0.33,\valtmp+0.66) -- (\tmp,\valtmp+0.99) -- (\tmp-0.33,\valtmp+0.99) -- (\tmp-0.33,\valtmp+0.66) -- cycle;
				  \pgfmathparse{50*rnd+50}
				  \edef\col{\pgfmathresult}
				  \fill[white!\col!black] (\tmp,\valtmp+0.33) -- (\tmp+0.33,\valtmp+0.33) -- (\tmp+0.33,\valtmp+0.66) -- (\tmp,\valtmp+0.99) -- (\tmp-0.33,\valtmp+0.99) -- (\tmp-0.33,\valtmp+0.66) -- cycle;	
			  }
		  }
		  \foreach \tmp in {0.99,1.98}{
			  \foreach \valtmp in {0.66,1.65}{
				  \draw[line width=\nonperiodic] (\tmp,\valtmp) -- (\tmp+0.33,\valtmp) -- (\tmp+0.33,\valtmp+0.33) -- (\tmp,\valtmp+0.66) -- (\tmp-0.33,\valtmp+0.66) -- (\tmp-0.33,\valtmp+0.33) -- cycle;
				  \pgfmathparse{50*rnd+50}
				  \edef\col{\pgfmathresult}
				  \fill[white!\col!black] (\tmp,\valtmp) -- (\tmp+0.33,\valtmp) -- (\tmp+0.33,\valtmp+0.33) -- 	(\tmp,\valtmp+0.66) -- (\tmp-0.33,\valtmp+0.66) -- (\tmp-0.33,\valtmp+0.33) -- cycle;
			  }
		  }
	\end{tikzpicture}

\end{minipage}
\begin{minipage}{0.48\textwidth}
	\centering
\begin{tabular}{rcccc}
		& & $\Nb{(X,i)}{\text{Ft}}$ & & \\
		& & \rotatebox{90}{$=$} & & \\
		\toprule
		$(X,i) \downarrow $  & $\Nb{(X,i)}{\text{V}}$ & $\Nb{(X,i)}{\text{Ed}}$ & \color{gray} $\Nb{(X,i)}{\text{Fa}}$ & $R^{(X,i)}_{\text{El}}$\\
		\midrule
		(V,1) & 13 & 15 & \color{gray}3 & 2 
		\\ 
		(Ft,1) = (Ed,1) & 10 & 11 & \color{gray}2 & 3 
		\\
		(El,1) = (Fa,1) & \color{gray}6 &  6 &  \color{gray}1 & \color{gray}1
		\\
		\bottomrule 
	\end{tabular}
\end{minipage}
\end{center}\vspace*{-0.5cm}
\caption{Sketch of structured periodic hexagonal mesh (left) and the relevant neighborhood topology numbers (right).}
	\label{fig:hexagons}
  \end{figure}  

For the periodic hexagon mesh, cf. \cref{fig:hexagons},
we again only have one type of vertices, edges and faces, respectively, i.e. $\{1\} = \mathcal{E}_0 = \mathcal{E}_1 = \mathcal{E}_2$. We obtain the dependency on $k$ in \ncdof and \nnze as shown in \cref{tab:hexagons:nnnpol}.

  \begin{table}[ht!]
	\begin{center}	
	\begin{NiceTabular}{rr@{~}r@{}r@{~}r@{}r|r@{~}r@{}r@{~}r@{}r@{~}r@{}r@{~}r@{}r}[colortbl-like]
\toprule
& \multicolumn{5}{c}{$\texttt{ncdof}/N_{El}$}& \multicolumn{9}{c}{$\texttt{nnze}/N_{El}$} \\
\midrule
DG$ $&$  \frac{1}{2} k^2$&$  + $&$  \frac{3}{2} k  $&$  + $&$ 1$&$  \frac{7}{4} k^4$&$  + $&$  \frac{21}{2} k^3$&$  + $&$  \frac{91}{4} k^2$&$  + $&$ 21 k  $&$  + $&$ 7$\\[0.5ex]\rowcolor{gray!20}\\[-2ex]
\rowcolor{gray!20}TDG2$ $&$ $&$ $&$ 2 k  $&$  + $&$ 1$&$ $&$ $&$ $&$ $&$ 28 k^2$&$  + $&$ 28 k  $&$  + $&$ 7$\\[0.5ex]\\[-2ex]
TDG1$ $&$ $&$ $&$  k  $&$  + $&$ 1$&$ $&$ $&$ $&$ $&$ 7 k^2$&$  + $&$ 14 k  $&$  + $&$ 7$\\[0.5ex]\rowcolor{gray!20}\\[-2ex]
\rowcolor{gray!20}HDG$ $&$ $&$ $&$ 3 k  $&$  + $&$ 3$&$ $&$ $&$ $&$ $&$ 33 k^2$&$  + $&$ 66 k  $&$  + $&$ 33$\\[0.5ex]\\[-2ex]
HHO$ $&$ $&$ $&$ 3 k  $&$ $&$ $&$ $&$ $&$ $&$ $&$ 33 k^2$&$ $&$ $&$ $&$ $\\[0.5ex]\rowcolor{gray!20}\\[-2ex]
\rowcolor{gray!20}VEM$ $&$ $&$ $&$ 3 k  $&$  - $&$ 1$&$ $&$ $&$ $&$ $&$ 33 k^2$&$  - $&$ 6 k  $&$  - $&$ 1$\\\bottomrule
\end{NiceTabular}

	\end{center}\vspace*{-0.5cm}
	\caption{Number of coupling unknowns and non-zero entries per element for different methods on a periodic quadrilateral mesh.} 
	\label{tab:hexagons:nnnpol}
	\end{table}

We observe that the scaling in $k$ turned -- compared to the triangle and quadrilateral case -- in favor of the Trefftz DG methods. There, the scaling with the $k^2$ term is lower now than for the skeleton-based discretizations HDG, HHO and VEM. In the lower order case, cf.  \cref{tab:hexagons:nnn}, we further observe that the HHO and VE methods still have slightly smaller numbers of non-zero entries in the linear system.

\begin{table}[!ht]
  \centering
  \begin{NiceTabular}{R{2.5cm}@{}R{1.05cm}@{}R{1.1cm}@{}R{1.15cm}@{}R{1.2cm}@{}R{1.25cm}@{}R{1.3cm}@{}R{1.35cm}@{}R{1.4cm}@{}R{1.45cm}@{}R{1.5cm}@{}}[colortbl-like]
    \toprule
\text{method}$~\downarrow$ $\setminus~k \rightarrow$ & \multicolumn{1}{r}{1}& \multicolumn{1}{r}{2}& \multicolumn{1}{r}{3}& \multicolumn{1}{r}{4}& \multicolumn{1}{r}{5}& \multicolumn{1}{r}{6}& \multicolumn{1}{r}{7}& \multicolumn{1}{r}{8}& \multicolumn{1}{r}{9}& \multicolumn{1}{r}{10}\\
   \midrule
$\texttt{ndof}_{\text{DG}}/N_{\text{el}}$&         3 &          6 &         10 &         15 &         21 &         28 &         36 &         45 &         55 &         66\\
\rowcolor{gray!20}$\texttt{ndof}_{\text{TDG2}}/N_{\text{el}}$&         3 &          5 &          7 &          9 &         11 &         13 &         15 &         17 &         19 &         21\\
\color{gray}$\texttt{ndof}_{\text{TDG1}}/N_{\text{el}}$&\color{gray}         2 & \color{gray}         3 & \color{gray}         4 & \color{gray}         5 & \color{gray}         6 & \color{gray}         7 & \color{gray}         8 & \color{gray}         9 & \color{gray}        10 & \color{gray}        11\\
\rowcolor{gray!20}$\texttt{ncdof}_{\text{HDG}}/N_{\text{el}}$&         6 &          9 &         12 &         15 &         18 &         21 &         24 &         27 &         30 &         33\\
\color{gray}$\texttt{ncdof}_{\text{HHO}}/N_{\text{el}}$&\color{gray}         3 & \color{gray}         6 & \color{gray}         9 & \color{gray}        12 & \color{gray}        15 & \color{gray}        18 & \color{gray}        21 & \color{gray}        24 & \color{gray}        27 & \color{gray}        30\\
\rowcolor{gray!20}$\texttt{ncdof}_{\text{VEM}}/N_{\text{el}}$&         2 &          5 &          8 &         11 &         14 &         17 &         20 &         23 &         26 &         29\\
    \bottomrule
\end{NiceTabular}

  \begin{NiceTabular}{R{2.5cm}@{}R{1.05cm}@{}R{1.1cm}@{}R{1.15cm}@{}R{1.2cm}@{}R{1.25cm}@{}R{1.3cm}@{}R{1.35cm}@{}R{1.4cm}@{}R{1.45cm}@{}R{1.5cm}@{}}[colortbl-like]
    \toprule
\text{method}$~\downarrow$ $\setminus~k \rightarrow$ & \multicolumn{1}{r}{1}& \multicolumn{1}{r}{2}& \multicolumn{1}{r}{3}& \multicolumn{1}{r}{4}& \multicolumn{1}{r}{5}& \multicolumn{1}{r}{6}& \multicolumn{1}{r}{7}& \multicolumn{1}{r}{8}& \multicolumn{1}{r}{9}& \multicolumn{1}{r}{10}\\
   \midrule
$\texttt{nnze}_{\text{DG}}/N_{\text{el}}$&        63 &        252 &        700 &       1575 &       3087 &       5488 &       9072 &      14175 &      21175 &      30492\\
\rowcolor{gray!20}$\texttt{nnze}_{\text{TDG2}}/N_{\text{el}}$&        63 &        175 &        343 &        567 &        847 &       1183 &       1575 &       2023 &       2527 &       3087\\
\color{gray}$\texttt{nnze}_{\text{TDG1}}/N_{\text{el}}$&\color{gray}        28 & \color{gray}        63 & \color{gray}       112 & \color{gray}       175 & \color{gray}       252 & \color{gray}       343 & \color{gray}       448 & \color{gray}       567 & \color{gray}       700 & \color{gray}       847\\
\rowcolor{gray!20}$\texttt{nnze}_{\text{HDG}}/N_{\text{el}}$&       132 &        297 &        528 &        825 &       1188 &       1617 &       2112 &       2673 &       3300 &       3993\\
\color{gray}$\texttt{nnze}_{\text{HHO}}/N_{\text{el}}$&\color{gray}        33 & \color{gray}       132 & \color{gray}       297 & \color{gray}       528 & \color{gray}       825 & \color{gray}      1188 & \color{gray}      1617 & \color{gray}      2112 & \color{gray}      2673 & \color{gray}      3300\\
\rowcolor{gray!20}$\texttt{nnze}_{\text{VEM}}/N_{\text{el}}$&        26 &        119 &        278 &        503 &        794 &       1151 &       1574 &       2063 &       2618 &       3239\\
    \bottomrule
\end{NiceTabular}

  \caption{\ncdof and \nnze per element for different methods on periodic hexagonal mesh for $k=1,\ldots,10$.}
  \label{tab:hexagons:nnn}
  \end{table}

\newpage
\subsection{3D: Tetrahedra}

\begin{figure}[ht!]
	\begin{center} 
		\begin{minipage}{0.48\textwidth}
 			\centering
	\includegraphics[width=0.5\textwidth]{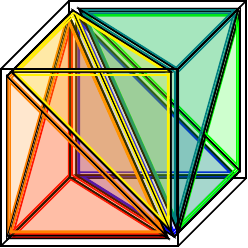}
\end{minipage}
\begin{minipage}{0.48\textwidth}
	\centering
\begin{tabular}{rcccc}
		& && $\Nb{(X,i)}{\text{Ft}}$ &  \\
		& && \rotatebox{90}{$=$} &  \\
		\toprule
		$(X,i) \downarrow $  & $\Nb{(X,i)}{\text{V}}$ & $\Nb{(X,i)}{\text{Ed}}$ & \color{gray} $\Nb{(X,i)}{\text{Fa}}$ & $R^{(X,i)}_{\text{El}}$\\
		\midrule
		(V,1) & 15 & 43 & 57 & $\frac16$ 
		\\[0.5ex] 
		(Ed,1) & 8 & 18 & 18 & $\frac23$
		\\[0.5ex] 
		(Ed,2) & 6 & 12 & 12 & $\frac12$
		\\[0.5ex] 
		(Ft,1) = (Fa,1) & 5 & 9 & 6 & 2
		\\[0.5ex] 
		(El,1) = (C,1) & \color{gray}4 &  \color{gray}6 & 4 & \color{gray}1
		\\
		\bottomrule 
	\end{tabular}
\end{minipage}	
\end{center}\vspace*{-0.5cm}
	\caption{Sketch of periodic unit cells based on Freudenthal decomposition of the unit cube into tetrahedra (left) and the relevant neighborhood topology numbers (right). Note that the diagonal edge in the decomposition has different neighborhood relations than all other edges.}
	\label{fig:tets}	
\end{figure} 
For the periodic tetrahedra mesh and a Freudenthal-decomposition of the unit cell, cf. \cref{fig:tets} we have only one type of vertices, faces and cells, respectively, i.e. $\{1\} = \mathcal{E}_0 = \mathcal{E}_2 = \mathcal{E}_3$, but two types of edges: 
The first type is shared by 6 elements, in \cref{fig:tets} these are the inner diagonal (from $(1,0,0)$ to $(0,1,1)$) and the outer edges of the cube ($\{(\xi,0,0)\}_{\xi\in[0,1]}$, $\{(0,\xi,0)\}_{\xi\in[0,1]}$, $\{(0,0,\xi)\}_{\xi\in[0,1]}$). 
The second type is shared by 4 elements and consists of the in-plane diagonal edges on the boundary of the unit cell. 
Thus $\{1,2\} = \mathcal{E}_1$. 
We obtain the dependency on $k$ in \ncdof and \nnze as shown in \cref{tab:tetra:nnnpol}. 

\begin{table}[ht!]
	\begin{center}	
	\begin{NiceTabular}{rr@{~}r@{}r@{~}r@{}r@{~}r@{}r|r@{~}r@{}r@{~}r@{}r@{~}r@{}r@{~}r@{}r@{~}r@{}r@{~}r@{}r}[colortbl-like]
\toprule
& \multicolumn{7}{c}{$\texttt{ncdof}/N_{El}$}& \multicolumn{13}{c}{$\texttt{nnze}/N_{El}$} \\
\midrule
DG$ $&$  \frac{1}{6} k^3$&$  + $&$  k^2$&$  + $&$  \frac{11}{6} k  $&$  + $&$ 1$&$  \frac{5}{36} k^6$&$  + $&$  \frac{5}{3} k^5$&$  + $&$  \frac{145}{18} k^4$&$  + $&$ 20 k^3$&$  + $&$  \frac{965}{36} k^2$&$  + $&$  \frac{55}{3} k  $&$  + $&$ 5$\\[0.5ex]\rowcolor{gray!20}\\[-2ex]
\rowcolor{gray!20}TDG2$ $&$ $&$ $&$  k^2$&$  + $&$ 2 k  $&$  + $&$ 1$&$ $&$ $&$ $&$ $&$ 5 k^4$&$  + $&$ 20 k^3$&$  + $&$ 30 k^2$&$  + $&$ 20 k  $&$  + $&$ 5$\\[0.5ex]\\[-2ex]
TDG1$ $&$ $&$ $&$  \frac{1}{2} k^2$&$  + $&$  \frac{3}{2} k  $&$  + $&$ 1$&$ $&$ $&$ $&$ $&$  \frac{5}{4} k^4$&$  + $&$  \frac{15}{2} k^3$&$  + $&$  \frac{65}{4} k^2$&$  + $&$ 15 k  $&$  + $&$ 5$\\[0.5ex]\rowcolor{gray!20}\\[-2ex]
\rowcolor{gray!20}HDG$ $&$ $&$ $&$  k^2$&$  + $&$ 3 k  $&$  + $&$ 2$&$ $&$ $&$ $&$ $&$  \frac{7}{2} k^4$&$  + $&$ 21 k^3$&$  + $&$  \frac{91}{2} k^2$&$  + $&$ 42 k  $&$  + $&$ 14$\\[0.5ex]\\[-2ex]
HHO$ $&$ $&$ $&$  k^2$&$  + $&$  k  $&$ $&$ $&$ $&$ $&$ $&$ $&$  \frac{7}{2} k^4$&$  + $&$ 7 k^3$&$  + $&$  \frac{7}{2} k^2$&$ $&$ $&$ $&$ $\\[0.5ex]\rowcolor{gray!20}\\[-2ex]
\rowcolor{gray!20}VEM$ $&$ $&$ $&$  k^2$&$  + $&$  \frac{1}{6} k  $&$  - $&$ 1$&$ $&$ $&$ $&$ $&$ 3 k^4$&$  + $&$ 12 k^3$&$  - $&$  \frac{21}{4} k^2$&$  - $&$  \frac{49}{4} k  $&$  + $&$ 5$\\\bottomrule
\end{NiceTabular}

	\end{center}\vspace*{-0.5cm}
	\caption{Number of coupling unknowns and non-zero entries per element for different methods on a periodic quadrilateral mesh.} 
	\label{tab:tetra:nnnpol} %
\end{table}

The number for the tetrahedron, the simplex in 3D, i.e. the geometry with the smallest amount of facets of an element and the smallest ratio between facets and elements in the mesh, is again in favor of the skeleton-based discretizations HDG, HHO and VEM. 
In \cref{tab:tetra:nnn} we list the concrete numbers of coupling dofs and non-zero entries in the linear system per element for the different methods on the periodic tetrahedra mesh for $k=1,\ldots,10$. The difference between the methods is less pronounced in 3D than it was in 2D and especially for the lower orders the Trefftz DG methods are competitive, with smaller numbers than HDG, but larger numbers than HHO and VEM. 

\begin{table}[!ht]
  \centering
  \begin{NiceTabular}{R{2.5cm}@{}R{1.05cm}@{}R{1.1cm}@{}R{1.15cm}@{}R{1.2cm}@{}R{1.25cm}@{}R{1.3cm}@{}R{1.35cm}@{}R{1.4cm}@{}R{1.45cm}@{}R{1.5cm}@{}}[colortbl-like]
    \toprule
\text{method}$~\downarrow$ $\setminus~k \rightarrow$ & \multicolumn{1}{r}{1}& \multicolumn{1}{r}{2}& \multicolumn{1}{r}{3}& \multicolumn{1}{r}{4}& \multicolumn{1}{r}{5}& \multicolumn{1}{r}{6}& \multicolumn{1}{r}{7}& \multicolumn{1}{r}{8}& \multicolumn{1}{r}{9}& \multicolumn{1}{r}{10}\\
   \midrule
$\texttt{ndof}_{\text{DG}}/N_{\text{el}}$&         4\hphantom{.0} &         10\hphantom{.0} &         20\hphantom{.0} &         35\hphantom{.0} &         56\hphantom{.0} &         84\hphantom{.0} &        120\hphantom{.0} &        165\hphantom{.0} &        220\hphantom{.0} &        286\hphantom{.0}\\
\rowcolor{gray!20}$\texttt{ndof}_{\text{TDG2}}/N_{\text{el}}$&         4\hphantom{.0} &          9\hphantom{.0} &         16\hphantom{.0} &         25\hphantom{.0} &         36\hphantom{.0} &         49\hphantom{.0} &         64\hphantom{.0} &         81\hphantom{.0} &        100\hphantom{.0} &        121\hphantom{.0}\\
\color{gray}$\texttt{ndof}_{\text{TDG1}}/N_{\text{el}}$&\color{gray}         3\hphantom{.0} & \color{gray}         6\hphantom{.0} & \color{gray}        10\hphantom{.0} & \color{gray}        15\hphantom{.0} & \color{gray}        21\hphantom{.0} & \color{gray}        28\hphantom{.0} & \color{gray}        36\hphantom{.0} & \color{gray}        45\hphantom{.0} & \color{gray}        55\hphantom{.0} & \color{gray}        66\hphantom{.0}\\
\rowcolor{gray!20}$\texttt{ncdof}_{\text{HDG}}/N_{\text{el}}$&         6\hphantom{.0} &         12\hphantom{.0} &         20\hphantom{.0} &         30\hphantom{.0} &         42\hphantom{.0} &         56\hphantom{.0} &         72\hphantom{.0} &         90\hphantom{.0} &        110\hphantom{.0} &        132\hphantom{.0}\\
\color{gray}$\texttt{ncdof}_{\text{HHO}}/N_{\text{el}}$&\color{gray}         2\hphantom{.0} & \color{gray}         6\hphantom{.0} & \color{gray}        12\hphantom{.0} & \color{gray}        20\hphantom{.0} & \color{gray}        30\hphantom{.0} & \color{gray}        42\hphantom{.0} & \color{gray}        56\hphantom{.0} & \color{gray}        72\hphantom{.0} & \color{gray}        90\hphantom{.0} & \color{gray}       110\hphantom{.0}\\
\rowcolor{gray!20}$\texttt{ncdof}_{\text{VEM}}/N_{\text{el}}$&       0.2 &        3.3 &        8.5 &       15.7 &       24.8 &         36\hphantom{.0} &       49.2 &       64.3 &       81.5 &      100.7\\
    \bottomrule
\end{NiceTabular}

  \begin{NiceTabular}{R{2.5cm}@{}R{1.05cm}@{}R{1.1cm}@{}R{1.15cm}@{}R{1.2cm}@{}R{1.25cm}@{}R{1.3cm}@{}R{1.35cm}@{}R{1.4cm}@{}R{1.45cm}@{}R{1.5cm}@{}}[colortbl-like]
    \toprule
\text{method}$~\downarrow$ $\setminus~k \rightarrow$ & \multicolumn{1}{r}{1}& \multicolumn{1}{r}{2}& \multicolumn{1}{r}{3}& \multicolumn{1}{r}{4}& \multicolumn{1}{r}{5}& \multicolumn{1}{r}{6}& \multicolumn{1}{r}{7}& \multicolumn{1}{r}{8}& \multicolumn{1}{r}{9}& \multicolumn{1}{r}{10}\\
   \midrule
$\texttt{nnze}_{\text{DG}}/N_{\text{el}}$&        80\hphantom{.0} &        500\hphantom{.0} &       2000\hphantom{.0} &       6125\hphantom{.0} &      15680\hphantom{.0} &      35280\hphantom{.0} &      72000\hphantom{.0} &     136125\hphantom{.0} &     242000\hphantom{.0} &     408980\hphantom{.0}\\
\rowcolor{gray!20}$\texttt{nnze}_{\text{TDG2}}/N_{\text{el}}$&        80\hphantom{.0} &        405\hphantom{.0} &       1280\hphantom{.0} &       3125\hphantom{.0} &       6480\hphantom{.0} &      12005\hphantom{.0} &      20480\hphantom{.0} &      32805\hphantom{.0} &      50000\hphantom{.0} &      73205\hphantom{.0}\\
\color{gray}$\texttt{nnze}_{\text{TDG1}}/N_{\text{el}}$&\color{gray}        45\hphantom{.0} & \color{gray}       180\hphantom{.0} & \color{gray}       500\hphantom{.0} & \color{gray}      1125\hphantom{.0} & \color{gray}      2205\hphantom{.0} & \color{gray}      3920\hphantom{.0} & \color{gray}      6480\hphantom{.0} & \color{gray}     10125\hphantom{.0} & \color{gray}     15125\hphantom{.0} & \color{gray}     21780\hphantom{.0}\\
\rowcolor{gray!20}$\texttt{nnze}_{\text{HDG}}/N_{\text{el}}$&       126\hphantom{.0} &        504\hphantom{.0} &       1400\hphantom{.0} &       3150\hphantom{.0} &       6174\hphantom{.0} &      10976\hphantom{.0} &      18144\hphantom{.0} &      28350\hphantom{.0} &      42350\hphantom{.0} &      60984\hphantom{.0}\\
\color{gray}$\texttt{nnze}_{\text{HHO}}/N_{\text{el}}$&\color{gray}        14\hphantom{.0} & \color{gray}       126\hphantom{.0} & \color{gray}       504\hphantom{.0} & \color{gray}      1400\hphantom{.0} & \color{gray}      3150\hphantom{.0} & \color{gray}      6174\hphantom{.0} & \color{gray}     10976\hphantom{.0} & \color{gray}     18144\hphantom{.0} & \color{gray}     28350\hphantom{.0} & \color{gray}     42350\hphantom{.0}\\
\rowcolor{gray!20}$\texttt{nnze}_{\text{VEM}}/N_{\text{el}}$&       2.5 &      103.5 &        488\hphantom{.0} &       1408\hphantom{.0} &     3187.5 &     6222.5 &      10981\hphantom{.0} &      18003\hphantom{.0} &    27900.5 &    41357.5\\
    \bottomrule
\end{NiceTabular}

  \caption{\ncdof and \nnze per element for different methods on periodic tetrahedra mesh for $k\!\!=\!\!1,\!..,\!10$ (rounded up to one decimal place).}  \label{tab:tetra:nnn}
\end{table}

\newpage  
\subsection{3D: Hexahedra}

\begin{figure}[ht!]
	\begin{center} 
		\begin{minipage}{0.48\textwidth}
			\centering
			
	\pgfmathsetmacro{\cubex}{3}
	\pgfmathsetmacro{\cubey}{3}
	\pgfmathsetmacro{\cubez}{-3}
	\begin{tikzpicture}
		\foreach \x in {0,1,2}{
			\foreach \y in {0,1,2}{
				\pgfmathparse{50*rnd+50}
				\edef\col{\pgfmathresult}
				\fill[white!\col!black, opacity=1] (\x*\cubex/3,\y*\cubey/3,0) -- (\x*\cubex/3+\cubex/3,\y*\cubey/3,0) -- (\x*\cubex/3+\cubex/3,\y*\cubey/3+\cubey/3,0) -- (\x*\cubex/3,\y*\cubey/3+\cubey/3,0) -- cycle;
				\pgfmathparse{50*rnd+50}
				\edef\col{\pgfmathresult}
				\fill[white!\col!black, opacity=1] (3*\cubex/3,\y*\cubey/3,\x*\cubez/3) -- (3*\cubex/3,\y*\cubey/3,\x*\cubez/3+\cubez/3) -- (3*\cubex/3,\y*\cubey/3+\cubey/3,\x*\cubez/3+\cubez/3) -- (3*\cubex/3,\y*\cubey/3+\cubey/3,\x*\cubez/3) -- cycle;
				\pgfmathparse{50*rnd+50}
				\edef\col{\pgfmathresult}
				\fill[white!\col!black, opacity=1] (\y*\cubex/3,3*\cubey/3,\x*\cubez/3) -- (\y*\cubex/3,3*\cubey/3,\x*\cubez/3+\cubez/3) -- (\y*\cubex/3+\cubex/3,3*\cubey/3,\x*\cubez/3+\cubez/3) -- (\y*\cubex/3+\cubex/3,3*\cubey/3,\x*\cubez/3) -- cycle;
			}
		}

		\foreach \x in {0,1,2,3}{
			\draw (\x*\cubex/3,0,0) -- (\x*\cubex/3,3*\cubey/3,0);
			\draw (0,\x*\cubey/3,0) -- (3*\cubex/3,\x*\cubey/3,0);
			\draw (\x*\cubex/3,3*\cubey/3,0) -- (\x*\cubex/3,3*\cubey/3,3*\cubez/3);
			\draw (3*\cubex/3,\x*\cubey/3,0) -- (3*\cubex/3,\x*\cubey/3,3*\cubez/3);
		}
		\draw (0,\cubey,\cubez) -- (\cubex,\cubey,\cubez) -- (\cubex,\cubey,\cubez) -- (\cubex,0,\cubez);
		\draw (0,\cubey,\cubez/3) -- (\cubex,\cubey,\cubez/3) -- (\cubex,0,\cubez/3);
		\draw (0,\cubey,2*\cubez/3) -- (\cubex,\cubey,2*\cubez/3) -- (\cubex,0,2*\cubez/3);
	\end{tikzpicture}

\end{minipage}
\begin{minipage}{0.48\textwidth}
	\centering
\begin{tabular}{rcccc}
		& && $\Nb{(X,i)}{\text{Ft}}$ &  \\
		& && \rotatebox{90}{$=$} &  \\
		\toprule
		$(X,i) \downarrow $  & $\Nb{(X,i)}{\text{V}}$ & $\Nb{(X,i)}{\text{Ed}}$ & \color{gray} $\Nb{(X,i)}{\text{Fa}}$ & $R^{(X,i)}_{\text{El}}$\\
		\midrule
		(V,1) & 27 & 54 & 12 & 1
		\\ 
		(Ed,1) & 18 & 33 & 20 & 3
		\\
		(Ft,1) = (Fa,1) & 12 & 20 & 11 & 3
		\\
		(El,1) = (C,1) & \color{gray}8 &  \color{gray}12 & 6 & \color{gray}1
		\\
		\bottomrule 
	\end{tabular}
\end{minipage}	
\end{center}\vspace*{-0.5cm}
	\caption{Sketch of periodic unit cells based on hexahedra (left) and the relevant neighborhood topology numbers (right).}
	\label{fig:hexa}	
\end{figure} 

For the periodic hexahedra mesh, cf. \cref{fig:hexa},
we again only have one type of vertices, edges, faces and cells respectively, i.e. $\{1\} = \mathcal{E}_0 = \mathcal{E}_1 = \mathcal{E}_2 = \mathcal{E}_3$. 
We obtain the dependency on $k$ in \ncdof and \nnze as shown in \cref{tab:hexa:nnnpol}.

\begin{table}[ht!]
	\begin{center}	
	\begin{NiceTabular}{rr@{~}r@{}r@{~}r@{}r@{~}r@{}r|r@{~}r@{}r@{~}r@{}r@{~}r@{}r@{~}r@{}r@{~}r@{}r@{~}r@{}r}[colortbl-like]
\toprule
& \multicolumn{7}{c}{$\texttt{ncdof}/N_{El}$}& \multicolumn{13}{c}{$\texttt{nnze}/N_{El}$} \\
\midrule
DG$ $&$  \frac{1}{6} k^3$&$  + $&$  k^2$&$  + $&$  \frac{11}{6} k  $&$  + $&$ 1$&$  \frac{7}{36} k^6$&$  + $&$  \frac{7}{3} k^5$&$  + $&$  \frac{203}{18} k^4$&$  + $&$ 28 k^3$&$  + $&$  \frac{1351}{36} k^2$&$  + $&$  \frac{77}{3} k  $&$  + $&$ 7$\\[0.5ex]\rowcolor{gray!20}\\[-2ex]
\rowcolor{gray!20}TDG2$ $&$ $&$ $&$  k^2$&$  + $&$ 2 k  $&$  + $&$ 1$&$ $&$ $&$ $&$ $&$ 7 k^4$&$  + $&$ 28 k^3$&$  + $&$ 42 k^2$&$  + $&$ 28 k  $&$  + $&$ 7$\\[0.5ex]\\[-2ex]
TDG1$ $&$ $&$ $&$  \frac{1}{2} k^2$&$  + $&$  \frac{3}{2} k  $&$  + $&$ 1$&$ $&$ $&$ $&$ $&$  \frac{7}{4} k^4$&$  + $&$  \frac{21}{2} k^3$&$  + $&$  \frac{91}{4} k^2$&$  + $&$ 21 k  $&$  + $&$ 7$\\[0.5ex]\rowcolor{gray!20}\\[-2ex]
\rowcolor{gray!20}HDG$ $&$ $&$ $&$  \frac{3}{2} k^2$&$  + $&$  \frac{9}{2} k  $&$  + $&$ 3$&$ $&$ $&$ $&$ $&$  \frac{33}{4} k^4$&$  + $&$  \frac{99}{2} k^3$&$  + $&$  \frac{429}{4} k^2$&$  + $&$ 99 k  $&$  + $&$ 33$\\[0.5ex]\\[-2ex]
HHO$ $&$ $&$ $&$  \frac{3}{2} k^2$&$  + $&$  \frac{3}{2} k  $&$ $&$ $&$ $&$ $&$ $&$ $&$  \frac{33}{4} k^4$&$  + $&$  \frac{33}{2} k^3$&$  + $&$  \frac{33}{4} k^2$&$ $&$ $&$ $&$ $\\[0.5ex]\rowcolor{gray!20}\\[-2ex]
\rowcolor{gray!20}VEM$ $&$ $&$ $&$  \frac{3}{2} k^2$&$  + $&$  \frac{3}{2} k  $&$  - $&$ 2$&$ $&$ $&$ $&$ $&$  \frac{33}{4} k^4$&$  + $&$  \frac{87}{2} k^3$&$  + $&$  \frac{45}{4} k^2$&$  - $&$ 54 k  $&$  + $&$ 18$\\\bottomrule
\end{NiceTabular}

	\end{center}\vspace*{-0.5cm}
	\caption{Number of coupling unknowns and non-zero entries per element for different methods on periodic hexahedra mesh.} 
	\label{tab:hexa:nnnpol}
	\end{table}

Similarly to the behavior in 2D, when going from the simplex to the tensor-product geometry, we increase the number of facets per element and the odds turn in favor of the cell-based discretizations, especially the Trefftz DG discretizations. In the higher order case $k \gg 1$, the Trefftz DG approaches outperform the skeleton-based discretizations in both the number of coupling dofs and the number of non-zero entries in the linear system. In the lower order case, and with regards to the HHO method also for all considered orders, cf.  \cref{tab:hexa:nnn}, we observe that the HHO and VE methods can still have slightly smaller numbers.

\begin{table}[!ht]
  \centering
  \begin{NiceTabular}{R{2.5cm}@{}R{1.05cm}@{}R{1.1cm}@{}R{1.15cm}@{}R{1.2cm}@{}R{1.25cm}@{}R{1.3cm}@{}R{1.35cm}@{}R{1.4cm}@{}R{1.45cm}@{}R{1.5cm}@{}}[colortbl-like]
    \toprule
\text{method}$~\downarrow$ $\setminus~k \rightarrow$ & \multicolumn{1}{r}{1}& \multicolumn{1}{r}{2}& \multicolumn{1}{r}{3}& \multicolumn{1}{r}{4}& \multicolumn{1}{r}{5}& \multicolumn{1}{r}{6}& \multicolumn{1}{r}{7}& \multicolumn{1}{r}{8}& \multicolumn{1}{r}{9}& \multicolumn{1}{r}{10}\\
   \midrule
$\texttt{ndof}_{\text{DG}}/N_{\text{el}}$&         4 &         10 &         20 &         35 &         56 &         84 &        120 &        165 &        220 &        286\\
\rowcolor{gray!20}$\texttt{ndof}_{\text{TDG2}}/N_{\text{el}}$&         4 &          9 &         16 &         25 &         36 &         49 &         64 &         81 &        100 &        121\\
\color{gray}$\texttt{ndof}_{\text{TDG1}}/N_{\text{el}}$&\color{gray}         3 & \color{gray}         6 & \color{gray}        10 & \color{gray}        15 & \color{gray}        21 & \color{gray}        28 & \color{gray}        36 & \color{gray}        45 & \color{gray}        55 & \color{gray}        66\\
\rowcolor{gray!20}$\texttt{ncdof}_{\text{HDG}}/N_{\text{el}}$&         9 &         18 &         30 &         45 &         63 &         84 &        108 &        135 &        165 &        198\\
\color{gray}$\texttt{ncdof}_{\text{HHO}}/N_{\text{el}}$&\color{gray}         3 & \color{gray}         9 & \color{gray}        18 & \color{gray}        30 & \color{gray}        45 & \color{gray}        63 & \color{gray}        84 & \color{gray}       108 & \color{gray}       135 & \color{gray}       165\\
\rowcolor{gray!20}$\texttt{ncdof}_{\text{VEM}}/N_{\text{el}}$&         1 &          7 &         16 &         28 &         43 &         61 &         82 &        106 &        133 &        163\\
    \bottomrule
\end{NiceTabular}

  \begin{NiceTabular}{R{2.5cm}@{}R{1.05cm}@{}R{1.1cm}@{}R{1.15cm}@{}R{1.2cm}@{}R{1.25cm}@{}R{1.3cm}@{}R{1.35cm}@{}R{1.4cm}@{}R{1.45cm}@{}R{1.5cm}@{}}[colortbl-like]
    \toprule
\text{method}$~\downarrow$ $\setminus~k \rightarrow$ & \multicolumn{1}{r}{1}& \multicolumn{1}{r}{2}& \multicolumn{1}{r}{3}& \multicolumn{1}{r}{4}& \multicolumn{1}{r}{5}& \multicolumn{1}{r}{6}& \multicolumn{1}{r}{7}& \multicolumn{1}{r}{8}& \multicolumn{1}{r}{9}& \multicolumn{1}{r}{10}\\
   \midrule
$\texttt{nnze}_{\text{DG}}/N_{\text{el}}$&       112 &        700 &       2800 &       8575 &      21952 &      49392 &     100800 &     190575 &     338800 &     572572\\
\rowcolor{gray!20}$\texttt{nnze}_{\text{TDG2}}/N_{\text{el}}$&       112 &        567 &       1792 &       4375 &       9072 &      16807 &      28672 &      45927 &      70000 &     102487\\
\color{gray}$\texttt{nnze}_{\text{TDG1}}/N_{\text{el}}$&\color{gray}        63 & \color{gray}       252 & \color{gray}       700 & \color{gray}      1575 & \color{gray}      3087 & \color{gray}      5488 & \color{gray}      9072 & \color{gray}     14175 & \color{gray}     21175 & \color{gray}     30492\\
\rowcolor{gray!20}$\texttt{nnze}_{\text{HDG}}/N_{\text{el}}$&       297 &       1188 &       3300 &       7425 &      14553 &      25872 &      42768 &      66825 &      99825 &     143748\\
\color{gray}$\texttt{nnze}_{\text{HHO}}/N_{\text{el}}$&\color{gray}        33 & \color{gray}       297 & \color{gray}      1188 & \color{gray}      3300 & \color{gray}      7425 & \color{gray}     14553 & \color{gray}     25872 & \color{gray}     42768 & \color{gray}     66825 & \color{gray}     99825\\
\rowcolor{gray!20}$\texttt{nnze}_{\text{VEM}}/N_{\text{el}}$&        27 &        435 &       1800 &       4878 &      10623 &      20187 &      34920 &      56370 &      86283 &     126603\\
    \bottomrule
\end{NiceTabular}

  \caption{\ncdof and \nnze per element for different methods on periodic hexahedra mesh for $k=1,\ldots,10$.}\label{tab:hexa:nnn}
  \end{table}

\newpage
\subsection{3D: Octahedron}

\begin{figure}[ht!]
	\begin{center} 
		\begin{minipage}{0.48\textwidth}
			\centering
	\includegraphics[width=0.5\textwidth]{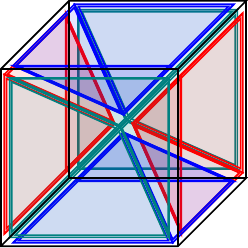}

\end{minipage}
\begin{minipage}{0.48\textwidth}
	\centering
\begin{tabular}{rcccc}
		& && $\Nb{(X,i)}{\text{Ft}}$ &  \\
		& && \rotatebox{90}{$=$} &  \\
		\toprule
		$(X,i) \downarrow $  & $\Nb{(X,i)}{\text{V}}$ & $\Nb{(X,i)}{\text{Ed}}$ & \color{gray} $\Nb{(X,i)}{\text{Fa}}$ & $R^{(X,i)}_{\text{El}}$\\
		\midrule
		(V,1) & 15 & 44 & 36 & $\frac13$ 
		\\[0.5ex] 
		(V,2) & 27 & 86 & 72 & $\frac23$ 
		\\[0.5ex] 
		(Ed,1) & 11 & 28 & 21 & $\frac83$
		\\[0.5ex] 
		(Ed,2) & 14 & 37 & 28 & 1
		\\[0.5ex] 
		(Ft,1) = (Fa,1) & 9 & 21 & 15 & 4
		\\[0.5ex] 
		(El,1) = (C,1) & \color{gray}6 &  \color{gray}12 & 8 & \color{gray}1
		\\
		\bottomrule 
	\end{tabular}
\end{minipage}	
\end{center}\vspace*{-0.5cm}
\caption{
 Sketch of periodic unit cells based on decomposition into three 
 octahedra
 (left) and the relevant neighborhood topology numbers (right).
 Note that always the two pyramids in the same color touch across the periodic boundary and form an octahedron.
}
	\label{fig:octa}	
\end{figure} 

Next, we increase the number of facets per element further by considering the periodic octahedra mesh, cf. \cref{fig:octa}. Here, we have only one type of faces and cells, respectively, i.e. $\{1\} = \mathcal{E}_2 = \mathcal{E}_3$, but two types of vertices and edges, i.e. $\{1,2\} = \mathcal{E}_0 = \mathcal{E}_1$. 
The different neighboring relations of the interior vertex and the interior edges, (V,1) and (Ed,1), and the outer edges and outer vertices, (V,2) and (Ed,2), respectively, are shown in \cref{fig:octa} and the accompanying table.
Outer and interior when referring to the vertices and edges are meant with respect to the unit cell.
We obtain the dependency on $k$ in \ncdof and \nnze as shown in
\cref{tab:octa:nnnpol}. 

\begin{table}[ht!]
	\begin{center}	
	\begin{NiceTabular}{rr@{~}r@{}r@{~}r@{}r@{~}r@{}r|r@{~}r@{}r@{~}r@{}r@{~}r@{}r@{~}r@{}r@{~}r@{}r@{~}r@{}r}[colortbl-like]
\toprule
& \multicolumn{7}{c}{$\texttt{ncdof}/N_{El}$}& \multicolumn{13}{c}{$\texttt{nnze}/N_{El}$} \\
\midrule
DG$ $&$  \frac{1}{6} k^3$&$  + $&$  k^2$&$  + $&$  \frac{11}{6} k  $&$  + $&$ 1$&$  \frac{1}{4} k^6$&$  + $&$ 3 k^5$&$  + $&$  \frac{29}{2} k^4$&$  + $&$ 36 k^3$&$  + $&$  \frac{193}{4} k^2$&$  + $&$ 33 k  $&$  + $&$ 9$\\[0.5ex]\rowcolor{gray!20}\\[-2ex]
\rowcolor{gray!20}TDG2$ $&$ $&$ $&$  k^2$&$  + $&$ 2 k  $&$  + $&$ 1$&$ $&$ $&$ $&$ $&$ 9 k^4$&$  + $&$ 36 k^3$&$  + $&$ 54 k^2$&$  + $&$ 36 k  $&$  + $&$ 9$\\[0.5ex]\\[-2ex]
TDG1$ $&$ $&$ $&$  \frac{1}{2} k^2$&$  + $&$  \frac{3}{2} k  $&$  + $&$ 1$&$ $&$ $&$ $&$ $&$  \frac{9}{4} k^4$&$  + $&$  \frac{27}{2} k^3$&$  + $&$  \frac{117}{4} k^2$&$  + $&$ 27 k  $&$  + $&$ 9$\\[0.5ex]\rowcolor{gray!20}\\[-2ex]
\rowcolor{gray!20}HDG$ $&$ $&$ $&$ 2 k^2$&$  + $&$ 6 k  $&$  + $&$ 4$&$ $&$ $&$ $&$ $&$ 15 k^4$&$  + $&$ 90 k^3$&$  + $&$ 195 k^2$&$  + $&$ 180 k  $&$  + $&$ 60$\\[0.5ex]\\[-2ex]
HHO$ $&$ $&$ $&$ 2 k^2$&$  + $&$ 2 k  $&$ $&$ $&$ $&$ $&$ $&$ $&$ 15 k^4$&$  + $&$ 30 k^3$&$  + $&$ 15 k^2$&$ $&$ $&$ $&$ $\\[0.5ex]\rowcolor{gray!20}\\[-2ex]
\rowcolor{gray!20}VEM$ $&$ $&$ $&$ 2 k^2$&$  + $&$  \frac{5}{3} k  $&$  - $&$  \frac{8}{3}$&$ $&$ $&$ $&$ $&$ 15 k^4$&$  + $&$ 54 k^3$&$  + $&$  \frac{20}{3} k^2$&$  - $&$ 72 k  $&$  + $&$  \frac{58}{3}$\\\bottomrule
\end{NiceTabular}

	\end{center}\vspace*{-0.5cm}
	\caption{Number of coupling unknowns and non-zero entries per element for different methods on periodic octahedra mesh.} 
	\label{tab:octa:nnnpol}
\end{table}

As we would expect from the tendency of the previous cases, the skeleton-based methods suffer from the increased number of coupled facets per cell rendering the Trefftz DG methods more efficient in almost all cases. Only for the orders 1--3 the Trefftz DG approach can be beaten by HHO (\ncdof and \nnze for order 1, \nnze also for order 2 and 3) and VEM (only order 1 and 2). 

\begin{table}[!ht]
  \centering

  \begin{NiceTabular}{R{2.5cm}@{}R{1.05cm}@{}R{1.1cm}@{}R{1.15cm}@{}R{1.2cm}@{}R{1.25cm}@{}R{1.3cm}@{}R{1.35cm}@{}R{1.4cm}@{}R{1.45cm}@{}R{1.5cm}@{}}[colortbl-like]
    \toprule
\text{method}$~\downarrow$ $\setminus~k \rightarrow$ & \multicolumn{1}{r}{1}& \multicolumn{1}{r}{2}& \multicolumn{1}{r}{3}& \multicolumn{1}{r}{4}& \multicolumn{1}{r}{5}& \multicolumn{1}{r}{6}& \multicolumn{1}{r}{7}& \multicolumn{1}{r}{8}& \multicolumn{1}{r}{9}& \multicolumn{1}{r}{10}\\
   \midrule
$\texttt{ndof}_{\text{DG}}/N_{\text{el}}$&         4\hphantom{.0} &         10\hphantom{.0} &         20\hphantom{.0} &         35\hphantom{.0} &         56\hphantom{.0} &         84\hphantom{.0} &        120\hphantom{.0} &        165\hphantom{.0} &        220\hphantom{.0} &        286\hphantom{.0}\\
\rowcolor{gray!20}$\texttt{ndof}_{\text{TDG2}}/N_{\text{el}}$&         4\hphantom{.0} &          9\hphantom{.0} &         16\hphantom{.0} &         25\hphantom{.0} &         36\hphantom{.0} &         49\hphantom{.0} &         64\hphantom{.0} &         81\hphantom{.0} &        100\hphantom{.0} &        121\hphantom{.0}\\
\color{gray}$\texttt{ndof}_{\text{TDG1}}/N_{\text{el}}$&\color{gray}         3\hphantom{.0} & \color{gray}         6\hphantom{.0} & \color{gray}        10\hphantom{.0} & \color{gray}        15\hphantom{.0} & \color{gray}        21\hphantom{.0} & \color{gray}        28\hphantom{.0} & \color{gray}        36\hphantom{.0} & \color{gray}        45\hphantom{.0} & \color{gray}        55\hphantom{.0} & \color{gray}        66\hphantom{.0}\\
\rowcolor{gray!20}$\texttt{ncdof}_{\text{HDG}}/N_{\text{el}}$&        12\hphantom{.0} &         24\hphantom{.0} &         40\hphantom{.0} &         60\hphantom{.0} &         84\hphantom{.0} &        112\hphantom{.0} &        144\hphantom{.0} &        180\hphantom{.0} &        220\hphantom{.0} &        264\hphantom{.0}\\
\color{gray}$\texttt{ncdof}_{\text{HHO}}/N_{\text{el}}$&\color{gray}         4\hphantom{.0} & \color{gray}        12\hphantom{.0} & \color{gray}        24\hphantom{.0} & \color{gray}        40\hphantom{.0} & \color{gray}        60\hphantom{.0} & \color{gray}        84\hphantom{.0} & \color{gray}       112\hphantom{.0} & \color{gray}       144\hphantom{.0} & \color{gray}       180\hphantom{.0} & \color{gray}       220\hphantom{.0}\\
\rowcolor{gray!20}$\texttt{ncdof}_{\text{VEM}}/N_{\text{el}}$&         1\hphantom{.0} &        8.7 &       20.3 &         36\hphantom{.0} &       55.7 &       79.3 &        107\hphantom{.0} &      138.7 &      174.3 &        214\hphantom{.0}\\
    \bottomrule
\end{NiceTabular}

  \begin{NiceTabular}{R{2.5cm}@{}R{1.05cm}@{}R{1.1cm}@{}R{1.15cm}@{}R{1.2cm}@{}R{1.25cm}@{}R{1.3cm}@{}R{1.35cm}@{}R{1.4cm}@{}R{1.45cm}@{}R{1.5cm}@{}}[colortbl-like]
    \toprule
\text{method}$~\downarrow$ $\setminus~k \rightarrow$ & \multicolumn{1}{r}{1}& \multicolumn{1}{r}{2}& \multicolumn{1}{r}{3}& \multicolumn{1}{r}{4}& \multicolumn{1}{r}{5}& \multicolumn{1}{r}{6}& \multicolumn{1}{r}{7}& \multicolumn{1}{r}{8}& \multicolumn{1}{r}{9}& \multicolumn{1}{r}{10}\\
   \midrule
$\texttt{nnze}_{\text{DG}}/N_{\text{el}}$&       144\hphantom{.0} &        900\hphantom{.0} &       3600\hphantom{.0} &      11025\hphantom{.0} &      28224\hphantom{.0} &      63504\hphantom{.0} &     129600\hphantom{.0} &     245025\hphantom{.0} &     435600\hphantom{.0} &     736164\hphantom{.0}\\
\rowcolor{gray!20}$\texttt{nnze}_{\text{TDG2}}/N_{\text{el}}$&       144\hphantom{.0} &        729\hphantom{.0} &       2304\hphantom{.0} &       5625\hphantom{.0} &      11664\hphantom{.0} &      21609\hphantom{.0} &      36864\hphantom{.0} &      59049\hphantom{.0} &      90000\hphantom{.0} &     131769\hphantom{.0}\\
\color{gray}$\texttt{nnze}_{\text{TDG1}}/N_{\text{el}}$&\color{gray}        81\hphantom{.0} & \color{gray}       324\hphantom{.0} & \color{gray}       900\hphantom{.0} & \color{gray}      2025\hphantom{.0} & \color{gray}      3969\hphantom{.0} & \color{gray}      7056\hphantom{.0} & \color{gray}     11664\hphantom{.0} & \color{gray}     18225\hphantom{.0} & \color{gray}     27225\hphantom{.0} & \color{gray}     39204\hphantom{.0}\\
\rowcolor{gray!20}$\texttt{nnze}_{\text{HDG}}/N_{\text{el}}$&       540\hphantom{.0} &       2160\hphantom{.0} &       6000\hphantom{.0} &      13500\hphantom{.0} &      26460\hphantom{.0} &      47040\hphantom{.0} &      77760\hphantom{.0} &     121500\hphantom{.0} &     181500\hphantom{.0} &     261360\hphantom{.0}\\
\color{gray}$\texttt{nnze}_{\text{HHO}}/N_{\text{el}}$&\color{gray}        60\hphantom{.0} & \color{gray}       540\hphantom{.0} & \color{gray}      2160\hphantom{.0} & \color{gray}      6000\hphantom{.0} & \color{gray}     13500\hphantom{.0} & \color{gray}     26460\hphantom{.0} & \color{gray}     47040\hphantom{.0} & \color{gray}     77760\hphantom{.0} & \color{gray}    121500\hphantom{.0} & \color{gray}    181500\hphantom{.0}\\
\rowcolor{gray!20}$\texttt{nnze}_{\text{VEM}}/N_{\text{el}}$&        23\hphantom{.0} &        574\hphantom{.0} &     2536.3 &       7134\hphantom{.0} &      15951\hphantom{.0} &    30931.3 &      54379\hphantom{.0} &      88958\hphantom{.0} &   137692.3 &     203966\hphantom{.0}\\
    \bottomrule
\end{NiceTabular}

  \caption{\ncdof and \nnze per element for different methods on periodic octahedra mesh for $k\!\!=\!\!1,\!..,\!10$ (rounded up to one decimal place).}
  \label{tab:octa:nnn}

  \end{table}

\newpage 
\subsection{3D: Truncated octahedron} \label{sec:truncatedoctahedron}

\begin{figure}[ht!]
	\begin{center} 
		\begin{minipage}{0.48\textwidth}
			\centering
	\pgfmathsetmacro{\cubex}{3}
	\pgfmathsetmacro{\cubey}{3}
	\pgfmathsetmacro{\cubez}{-3}
    \resizebox{0.45\textwidth}{!}{
    \begin{tikzpicture}[tdplot_main_coords,line cap=round,line join=round]
    \path foreach \Y in {0,1,2} {foreach \X in {0,1,2} 
      {({2*\Y}, {2*\X}, {-sqrt(2)*\X-sqrt(2)*\Y})
      node{\pgfmathtruncatemacro{\Z}{\X+\Y}
      \ifodd\Z
      \usebox{\TruncatedOctahedronOrange}
      \else
      \usebox{\TruncatedOctahedronTeal}
      \fi} }};
    \end{tikzpicture}}

\end{minipage}
\begin{minipage}{0.48\textwidth}
	\centering
\begin{tabular}{rcccc}
		& && $\Nb{(X,i)}{\text{Ft}}$ &  \\
		& && \rotatebox{90}{$=$} &  \\
		\toprule
		$(X,i) \downarrow $  & $\Nb{(X,i)}{\text{V}}$ & $\Nb{(X,i)}{\text{Ed}}$ & \color{gray} $\Nb{(X,i)}{\text{Fa}}$ & $R^{(X,i)}_{\text{El}}$\\
		\midrule
		(V,1) & 71 & 116 & 50 & 5
		\\
		(Ed,1) & 58 & 93 & 39 & 12
		\\
		(Ft,1) = (Fa,1) & 42 & 66 & 27 & 3
		\\
		(Ft,2) = (Fa,2) & 44 & 68 & 27 & 4
		\\
		(El,1) = (C,1) & \color{gray}24 &  \color{gray}36 & 14 & \color{gray}1
		\\
		\bottomrule 
	\end{tabular}
\end{minipage}	
\end{center}\vspace*{-0.5cm}	
	\caption{Sketch of periodic filling based on truncated octahedra (left) and the relevant neighborhood topology numbers (right).}
	\label{fig:truncoct}	
\end{figure} 
Finally, we turn our attention to the most extreme case in terms of the number of facets of each cell, the periodic truncated octahedron mesh, cf. \cref{fig:truncoct}. Note that these periodic structures can occur for example as the Voronoi tessellation of a body-centered cubic lattice crystal structure.
Here, we only have one type of vertices, edges, and cells respectively, i.e. $\{1\} = \mathcal{E}_0  = \mathcal{E}_1 = \mathcal{E}_3$. But two types of faces, $\{1,2\} = \mathcal{E}_2$, the hexagonal and the square faces, respectively. We obtain the dependency on $k$ in \ncdof and \nnze as shown in \cref{tab:truncoct:nnnpol}.
\begin{table}[ht!]
	\begin{center}	
	\begin{NiceTabular}{rr@{~}r@{}r@{~}r@{}r@{~}r@{}r|r@{~}r@{}r@{~}r@{}r@{~}r@{}r@{~}r@{}r@{~}r@{}r@{~}r@{}r}[colortbl-like]
\toprule
& \multicolumn{7}{c}{$\texttt{ncdof}/N_{El}$}& \multicolumn{13}{c}{$\texttt{nnze}/N_{El}$} \\
\midrule
DG$ $&$  \frac{1}{6} k^3$&$  + $&$  k^2$&$  + $&$  \frac{11}{6} k  $&$  + $&$ 1$&$  \frac{5}{12} k^6$&$  + $&$ 5 k^5$&$  + $&$  \frac{145}{6} k^4$&$  + $&$ 60 k^3$&$  + $&$  \frac{965}{12} k^2$&$  + $&$ 55 k  $&$  + $&$ 15$\\[0.5ex]\rowcolor{gray!20}\\[-2ex]
\rowcolor{gray!20}TDG2$ $&$ $&$ $&$  k^2$&$  + $&$ 2 k  $&$  + $&$ 1$&$ $&$ $&$ $&$ $&$ 15 k^4$&$  + $&$ 60 k^3$&$  + $&$ 90 k^2$&$  + $&$ 60 k  $&$  + $&$ 15$\\[0.5ex]\\[-2ex]
TDG1$ $&$ $&$ $&$  \frac{1}{2} k^2$&$  + $&$  \frac{3}{2} k  $&$  + $&$ 1$&$ $&$ $&$ $&$ $&$  \frac{15}{4} k^4$&$  + $&$  \frac{45}{2} k^3$&$  + $&$  \frac{195}{4} k^2$&$  + $&$ 45 k  $&$  + $&$ 15$\\[0.5ex]\rowcolor{gray!20}\\[-2ex]
\rowcolor{gray!20}HDG$ $&$ $&$ $&$  \frac{7}{2} k^2$&$  + $&$  \frac{21}{2} k  $&$  + $&$ 7$&$ $&$ $&$ $&$ $&$  \frac{189}{4} k^4$&$  + $&$  \frac{567}{2} k^3$&$  + $&$  \frac{2457}{4} k^2$&$  + $&$ 567 k  $&$  + $&$ 189$\\[0.5ex]\\[-2ex]
HHO$ $&$ $&$ $&$  \frac{7}{2} k^2$&$  + $&$  \frac{7}{2} k  $&$ $&$ $&$ $&$ $&$ $&$ $&$  \frac{189}{4} k^4$&$  + $&$  \frac{189}{2} k^3$&$  + $&$  \frac{189}{4} k^2$&$ $&$ $&$ $&$ $\\[0.5ex]\rowcolor{gray!20}\\[-2ex]
\rowcolor{gray!20}VEM$ $&$ $&$ $&$  \frac{7}{2} k^2$&$  + $&$  \frac{17}{2} k  $&$  - $&$ 6$&$ $&$ $&$ $&$ $&$  \frac{189}{4} k^4$&$  + $&$  \frac{749}{2} k^3$&$  + $&$  \frac{2105}{4} k^2$&$  - $&$ 672 k  $&$  + $&$ 150$\\\bottomrule
\end{NiceTabular}

	\end{center}\vspace*{-0.5cm}
	\caption{Number of coupling unknowns and non-zero entries per element for different methods on periodic truncated octahedra mesh.} 
	\label{tab:truncoct:nnnpol}
\end{table}

The number of lower-dimensional entities (faces, edges and vertices) per cell is now so large that the Trefftz DG methods are clearly superior to the skeleton-based methods, cf. \cref{tab:truncoct:nnn}. Only in the case $k=1$ the HHO method has a slightly smaller number of non-zero entries in the linear system. In all other cases, the Trefftz DG methods have the (much) smaller numbers.
\begin{table}[!ht]
  \centering
  \begin{NiceTabular}{R{2.5cm}@{}R{1.05cm}@{}R{1.1cm}@{}R{1.15cm}@{}R{1.2cm}@{}R{1.25cm}@{}R{1.3cm}@{}R{1.35cm}@{}R{1.4cm}@{}R{1.45cm}@{}R{1.5cm}@{}}[colortbl-like]
    \toprule
\text{method}$~\downarrow$ $\setminus~k \rightarrow$ & \multicolumn{1}{r}{1}& \multicolumn{1}{r}{2}& \multicolumn{1}{r}{3}& \multicolumn{1}{r}{4}& \multicolumn{1}{r}{5}& \multicolumn{1}{r}{6}& \multicolumn{1}{r}{7}& \multicolumn{1}{r}{8}& \multicolumn{1}{r}{9}& \multicolumn{1}{r}{10}\\
   \midrule
$\texttt{ndof}_{\text{DG}}/N_{\text{el}}$&         4 &         10 &         20 &         35 &         56 &         84 &        120 &        165 &        220 &        286\\
\rowcolor{gray!20}$\texttt{ndof}_{\text{TDG2}}/N_{\text{el}}$&         4 &          9 &         16 &         25 &         36 &         49 &         64 &         81 &        100 &        121\\
\color{gray}$\texttt{ndof}_{\text{TDG1}}/N_{\text{el}}$&\color{gray}         3 & \color{gray}         6 & \color{gray}        10 & \color{gray}        15 & \color{gray}        21 & \color{gray}        28 & \color{gray}        36 & \color{gray}        45 & \color{gray}        55 & \color{gray}        66\\
\rowcolor{gray!20}$\texttt{ncdof}_{\text{HDG}}/N_{\text{el}}$&        21 &         42 &         70 &        105 &        147 &        196 &        252 &        315 &        385 &        462\\
\color{gray}$\texttt{ncdof}_{\text{HHO}}/N_{\text{el}}$&\color{gray}         7 & \color{gray}        21 & \color{gray}        42 & \color{gray}        70 & \color{gray}       105 & \color{gray}       147 & \color{gray}       196 & \color{gray}       252 & \color{gray}       315 & \color{gray}       385\\
\rowcolor{gray!20}$\texttt{ncdof}_{\text{VEM}}/N_{\text{el}}$&         6 &         25 &         51 &         84 &        124 &        171 &        225 &        286 &        354 &        429\\
    \bottomrule
\end{NiceTabular}

  \begin{NiceTabular}{R{2.5cm}@{}R{1.05cm}@{}R{1.1cm}@{}R{1.15cm}@{}R{1.2cm}@{}R{1.25cm}@{}R{1.3cm}@{}R{1.35cm}@{}R{1.4cm}@{}R{1.45cm}@{}R{1.5cm}@{}}[colortbl-like]
    \toprule
\text{method}$~\downarrow$ $\setminus~k \rightarrow$ & \multicolumn{1}{r}{1}& \multicolumn{1}{r}{2}& \multicolumn{1}{r}{3}& \multicolumn{1}{r}{4}& \multicolumn{1}{r}{5}& \multicolumn{1}{r}{6}& \multicolumn{1}{r}{7}& \multicolumn{1}{r}{8}& \multicolumn{1}{r}{9}& \multicolumn{1}{r}{10}\\
   \midrule
$\texttt{nnze}_{\text{DG}}/N_{\text{el}}$&       240 &       1500 &       6000 &      18375 &      47040 &     105840 &     216000 &     408375 &     726000 &    1226940\\
\rowcolor{gray!20}$\texttt{nnze}_{\text{TDG2}}/N_{\text{el}}$&       240 &       1215 &       3840 &       9375 &      19440 &      36015 &      61440 &      98415 &     150000 &     219615\\
\color{gray}$\texttt{nnze}_{\text{TDG1}}/N_{\text{el}}$&\color{gray}       135 & \color{gray}       540 & \color{gray}      1500 & \color{gray}      3375 & \color{gray}      6615 & \color{gray}     11760 & \color{gray}     19440 & \color{gray}     30375 & \color{gray}     45375 & \color{gray}     65340\\
\rowcolor{gray!20}$\texttt{nnze}_{\text{HDG}}/N_{\text{el}}$&      1701 &       6804 &      18900 &      42525 &      83349 &     148176 &     244944 &     382725 &     571725 &     823284\\
\color{gray}$\texttt{nnze}_{\text{HHO}}/N_{\text{el}}$&\color{gray}       189 & \color{gray}      1701 & \color{gray}      6804 & \color{gray}     18900 & \color{gray}     42525 & \color{gray}     83349 & \color{gray}    148176 & \color{gray}    244944 & \color{gray}    382725 & \color{gray}    571725\\
\rowcolor{gray!20}$\texttt{nnze}_{\text{VEM}}/N_{\text{el}}$&       426 &       4663 &      16809 &      41946 &      86290 &     157191 &     263133 &     413734 &     619746 &     893055\\
    \bottomrule
\end{NiceTabular}

  \caption{\ncdof and \nnze per element for different methods on periodic truncated octahedron mesh for $k=1,\ldots,10$.} \label{tab:truncoct:nnn}
  \end{table}
\newpage

\section{Conclusion} \label{sec:conclusion}
Polytopal finite elements have become increasingly popular in the last years. While traditional $H^1$-conforming finite element methods are not directly generalized from simple geometries like simplices or cuboids, Discontinuous Galerkin methods, their hybrid versions and the virtual element methods can be applied on general polytopal meshes. 
In this work, we put the question aside which methods can be applied on polytopal meshes with a completely established theory for which class of PDE problems and instead focus on the question of efficiency -- under the assumption that the methods can be applied (equally well) and that the problem at hand is a scalar one.

We considered different types of periodic meshes in 2D and 3D ranging from simplex meshes to polytopal meshes with many more (faces,) edges and vertices in the mesh than elements. For these types of meshes we investigated the number of coupling unknowns and non-zero entries in the linear system for different methods. As one could expect, the methods' performances change depending on the ratio between geometrical entities of different dimensions. 
Hybrid Discontinuous Galerkin methods, Hybrid High-Order methods and Virtual Element methods are based on the idea of moving all coupling unknowns to the mesh skeleton. This gains efficiency over the classical element-based Discontinuous Galerkin methods for simple geometries like simplices or cuboids, where the number of lower dimensional entities is small or moderate compared to the number of elements. Especially on simplices, this leads to a large reduction in the number of coupling unknowns and non-zero entries in the linear system compared to standard Discontinuous Galerkin methods and also slightly outperforms the Trefftz Discontinuous Galerkin methods (for second order PDEs).
However, in the context of polytopal meshes and geometrical configurations with more facets (and lower dimensional entities) in the mesh as they appear for instance in dual-meshes or Voronoi-type constructions, 
 as seen in the periodic octahedron or the periodic truncated octahedron, are more relevant. In these configurations, we observe that the element-based Trefftz Discontinuous Galerkin methods outperform the skeleton-based methods, especially for higher orders.

\section*{Data availability statement}

The tables with the numerical results are generated by a simple script that we make available, cf. \url{https://doi.org/10.25625/H3AN3A} (see \cite{data}).

\section*{Acknowledgements}
This research was funded in part by the Austrian Science Fund (FWF) project \href{https://doi.org/10.55776/F65}{10.55776/F65}.

\bibliographystyle{plain}
\bibliography{LSZ}
\end{document}